\documentclass[12pt]{amsart}
\usepackage{graphicx, overpic, epsf,amssymb,amscd,times}

\newtheorem{thm}{Theorem}[section]

\newtheorem{prop}[thm]{Proposition}

\newtheorem{lemma}[thm]{Lemma}
\newtheorem{cor}[thm]{Corollary}

\newtheorem{q}[thm]{Question}
\newtheorem{rmk}[thm]{Remark}

\newcommand{\R}{\mathbb{R}}
\newcommand{\Z}{\mathbb{Z}}

\renewcommand{\P}{\mathbb{P}}
\renewcommand{\H}{\mathbb{H}}
\newcommand{\bdry}{\partial}

\newcommand{\s}{\vskip.1in}
\newcommand{\n}{\noindent}

\newcommand{\be}{\begin{enumerate}}
\newcommand{\ee}{\end{enumerate}}
\newcommand{\veer}{Veer(S,\bdry S)}
\newcommand{\dehn}{Dehn^+(S,\bdry S)}

\newcommand{\rot}{rot}

\numberwithin{equation}{subsection}

\begin{document}

\title{Right-veering diffeomorphisms of compact surfaces with boundary II}

\author{Ko Honda}
\address{University of Southern California, Los Angeles, CA 90089}
\email{khonda@math.usc.edu}
\urladdr{http://rcf.usc.edu/\char126 khonda}

\author{William H. Kazez}
\address{University of Georgia, Athens, GA 30602} \email{will@math.uga.edu}
\urladdr{http://www.math.uga.edu/\char126 will}

\author{Gordana Mati\'c}
\address{University of Georgia, Athens, GA 30602} \email{gordana@math.uga.edu}
\urladdr{http://www.math.uga.edu/\char126 gordana}

\date{This version: April 20, 2008. (The pictures are in color.)}



\thanks{KH supported by an Alfred P.\ Sloan Fellowship and an NSF
CAREER Award (DMS-0237386); GM supported by NSF grant DMS-0410066;
WHK supported by NSF grant DMS-0406158.}

\begin{abstract}
We continue our study of the monoid of right-veering diffeomorphisms
on a compact oriented surface with nonempty boundary, introduced in
\cite{HKM2}.  We conduct a detailed study of the case when the
surface is a punctured torus; in particular, we exhibit the
difference between the monoid of right-veering diffeomorphisms and
the monoid of products of positive Dehn twists, with the help of the
{\em Rademacher function}.  We then generalize to the braid group
$B_n$ on $n$ strands by relating the {\em signature} and the {\em
Maslov index}.  Finally, we discuss the symplectic fillability in
the pseudo-Anosov case by comparing with the work of
Roberts~\cite{Ro1,Ro2}.
\end{abstract}

\maketitle

\section{Introduction}

The paper \cite{HKM2} introduced the study of right-veering
diffeomorphisms on a compact oriented surface with nonempty boundary
(sometimes called a ``bordered surface'').  This paper continues the
investigations initiated in \cite{HKM2}.

Let $Aut(S,\bdry S)$ be the isotopy classes of diffeomorphisms of a
bordered surface $S$ which restrict to the identity on the boundary,
$\veer$ be the monoid of right-veering diffeomorphisms of $S$, and
$\dehn$ be the monoid of products of positive Dehn twists.  (In
particular, $id$ is in both.) Recall that, by the work of
Giroux~\cite{Gi2}, there is a 1--1 correspondence between
isomorphism classes of open book decompositions modulo stabilization
and isomorphism classes of contact structures on closed
$3$-manifolds. (Open books were introduced into contact geometry
much earlier by Thurston and Winkelnkemper~\cite{TW}.) If $h\in
Aut(S,\bdry S)$, let us write $(S,h)$ to denote, by slight abuse of
notation, either the open book decomposition or the corresponding
adapted contact structure.   The main result of \cite{HKM2} is that
a contact 3-manifold $(M,\xi)$ is tight if and only if all its
adapted open book decompositions have right-veering monodromy. Here
$M$ is closed and oriented, and $\xi$ is cooriented. On the other
hand, Giroux~\cite{Gi2} showed that $(M,\xi)$ is Stein fillable if
and only if there is an adapted open book decomposition with
monodromy $h\in \dehn$.  In order to understand the difference
between tight and Stein fillable contact structures, as well as the
symplectically fillable contact structures, which sit in between the
two, we need to understand the difference between $\veer$ and
$\dehn$.

One of the goals of this paper is to give an analysis of the
difference between $\veer$ and $\dehn$ for the once-punctured torus
$S$. The {\em Rademacher function} $\Phi$ and the {\em rotation
number} $rot$, defined in Sections~\ref{rade} and~\ref{link}, taken
together, are effective at distinguishing large swathes of $\veer$
that are not in $\dehn$.  Our first theorem is the following:

\begin{thm}\label{thm: veer not in dehn}
Let $S$ be a once-punctured torus and $h\in Aut(S,\bdry S)$. If
$rot(h)\geq {1\over 2}$ and $-\Phi(h)\geq 10~rot(h)$, then $h$ is in
$\veer$ but not in $\dehn$.
\end{thm}

\begin{proof}
This follows from Lemma~\ref{rotationlarge} and
Theorem~\ref{torusmain}.  (The lemma and the theorem are stated in
terms of $\sigma$ in $B_3$, the braid group $B_3$ on $3$ strands.
See below for the discussion of $B_3\cong Aut(S,\bdry S)$.)
\end{proof}

Theorem~\ref{torusmain} is, to a large extent, a consequence of the
fact that the linking number is positive on nontrivial elements of
$\dehn$.  However, we also give evidence that the linking number is
only a ``first-order'' invariant, in the sense that there are
elements in $\veer-\dehn$ which cannot be measured by this
technique, and require finer analysis.  Examples of this are given
in Section~\ref{examples}.

We will also generalize Theorem~\ref{torusmain} to the case of the
braid group $B_n$ on $n$ strands.  If $S$ is a double branched cover
of the disk, branched at $n$ points, then the {\em hyperelliptic}
mapping class group $HypAut(S,\bdry S)$ is the subgroup of
$Aut(S,\bdry S)$, consisting of diffeomorphisms that commute with
the hyperelliptic involution.  Equivalently, it is the image of
$B_n$ in $Aut(S,\bdry S)$.  When $n=3$, $S$ is a punctured torus and
half-twists about arcs connecting branch points lift to Dehn twists
that generate $Aut(S,\bdry S)$.  Hence $B_3$ can be identified with
$HypAut(S,\bdry S)=Aut(S,\bdry S)$.  On the other hand, for $n>3$,
$HypAut(S,\bdry S)$ is a proper subgroup of $Aut(S,\bdry S)$.

Now, an element $\sigma\in B_n$ which is a product of conjugates of
the standard positive half-twist generators is said to be {\em
quasipositive}. The monoid of quasipositive braids corresponds to
the monoid of products of positive Dehn twists, each of which is in
$HypAut(S,\bdry S)$. Observe that the monoid of quasipositive braids
strictly contains the monoid of {\em positive} braids, i.e., those
which are positive products of the standard generators
$\sigma_1,\dots,\sigma_{n-1}$ of $B_n$.  In Section
~\ref{braidgroup} we prove Theorem~\ref{sign+maslov}, following the
works of Gambaudo-Ghys \cite{GG1,GG2}. After lifting the action of
$HypAut(S,\bdry S)$ on homology to $\widetilde{Sp}(2n,\R)$, this
theorem describes a relationship between the {\em signature} of the
braid closure and the {\em Maslov index} of a corresponding ``lift''
to $\widetilde{Sp}(2n,\R)$. Corollary~\ref{cor: restriction} is then
an incarnation of the fact that the linking number is positive on
nontrivial quasipositive braids.

We then focus our attention to the question of which right-veering
monodromy maps $h$ correspond to tight contact structures.  In the
pseudo-Anosov case we have the following result, which is proved in
Section~\ref{section: weak fillability}:

\begin{thm}  \label{cgeq1}
Let $S$ be a bordered surface with connected boundary and $h$ be
pseudo-Anosov with fractional Dehn twist coefficient $c$.  If $c\geq
1$, then $(S,h)$ is a isotopic to a perturbation of a taut
foliation.  Hence $(S,h)$ is (weakly) symplectically fillable and
universally tight if $c\geq 1$.
\end{thm}

Hatcher~\cite{Ha} and Roberts~\cite{Ro1, Ro2} constructed
non-finite-depth taut foliations on certain Dehn fillings of
punctured surface bundles. (Hatcher's work was for punctured torus
bundles, which in turn was generalized by Roberts to all punctured
surface bundles with one boundary puncture.) Theorem~\ref{cgeq1}
follows from showing that the contact structure $(S,h)$ adapted to
the open book is isotopic to perturbations of the Hatcher-Roberts
taut foliations, using techniques developed in \cite{HKM1}.

We are now left to analyze $(S,h)$ when $h$ is pseudo-Anosov and the
fractional Dehn twist coefficient satisfies $0<c<1$.  (Recall that
if $c\leq 0$ then $(S,h)$ is overtwisted by Proposition~3.1 of
\cite{HKM2}.)  For example, when $S$ is a punctured torus, we are
concerned with $c={1\over 2}$.  In the paper~\cite{HKM3}, we prove,
using Heegaard Floer homology, that $(S,h)$ is tight if $c={1\over
2}$. This shows that if $S$ is a punctured torus and $h$ is
pseudo-Anosov, then $(S,h)$ is tight if and only if $h\in \veer$.

\section{$\veer$ vs. $\dehn$ on the punctured torus}

In this section we explain how to exhibit $h\in \veer$ that are not
products of positive Dehn twists, primarily via a combination of the
{\em Rademacher function} and the {\em rotation number}.

\subsection{Preliminaries}

We discuss some preliminary notions, partly to fix terminology.

Let $S$ be the once-punctured torus and $T$ be the torus.  There is
a short exact sequence
\begin{equation}\label{exact}
0\to {\Z} \to Aut(S, \partial S) \to Aut(T) \to 1,
\end{equation}
where the generator of $\Z$ is mapped to a positive Dehn twist
$R_{\bdry S}$ about $\partial S$. (In general, we use the notation
$R_\gamma$ to denote a positive Dehn twist about a closed curve
$\gamma$.) The group $Aut(T)$ is isomorphic to $SL(2,\Z)$, and is
generated by
$$A=\begin{pmatrix} 0 & 1\\ -1 & 0\end{pmatrix}, ~~~ B=
\begin{pmatrix} 1 & -1 \\ 1& 0 \end{pmatrix}.$$
Now, $Aut(S,\bdry S)$ can be identified with the (Artin) braid group
$B_3$ on 3 strands. Denote the generators of $B_3$ by $\sigma_1$ and
$\sigma_2$, corresponding to positive half-twists about strands 1, 2
and strands 2, 3.  We then have the relation
$\sigma_1\sigma_2\sigma_1 = \sigma_2\sigma_1\sigma_2$. If we view
the punctured torus as a 2-fold branched cover of the disk with 3
branch points, then the positive half-twists on $B_3$ lift to
positive Dehn twists on the punctured torus. More precisely, we
choose the images $\overline\sigma_i$ of $\sigma_i$ in $SL(2,\Z)$ to
be
$$\overline\sigma_1=\begin{pmatrix} 1 & 0 \\ -1 & 1\end{pmatrix},
~~~\overline\sigma_2 =\begin{pmatrix} 1 & 1\\ 0 & 1\end{pmatrix}.$$
Since $\overline\sigma_1 \overline\sigma_2 \overline\sigma_1
=\overline\sigma_2 \overline\sigma_1 \overline\sigma_2$ in
$SL(2,\Z)$, the map $\sigma_i \mapsto \overline\sigma_i$ induces a
homomorphism $B_3\to SL(2,\Z)$.  We have
$A=\overline\sigma_1\overline\sigma_2\overline\sigma_1$ and $B=
\overline\sigma_1^{-1}\overline\sigma_2^{-1}$.  In the short exact
sequence \ref{exact}, $R_{\bdry S} \mapsto
(\sigma_1\sigma_2\sigma_1)^4 \mapsto A^4$, which is the identity
matrix in $SL(2,\Z)$.

Elements of $SL(2, \Z)$ are grouped into three categories:
reducible, periodic and Anosov.  We will interpret the results from
\cite{HKM2} to determine which elements $h$ of $Aut(S,
\partial S)$ are right-veering. By a slight abuse of terminology, we
will often say ``$h$ is Anosov'' to mean ``$\overline{h}$ is an
Anosov diffeomorphism''.

If $h$ is Anosov, then $\overline{h}$ has a pair of irrational
eigenvalues $\lambda_1,\lambda_2$ that are both positive or both
negative. In either case, there are two prongs of the stable
lamination. If the $\lambda_i$ are positive, then the fractional
Dehn twist coefficient $c$ is an integer $n$ and the prongs are
fixed; if the $\lambda_i$ are negative, then $c$ is a half-integer
$n+{1\over 2}$ and the prongs are switched.  According to
\cite{HKM2}, an Anosov diffeomorphism $h$ is right-veering if and
only if $c\geq {1\over 2}$.

If $h$ is periodic, then $h$ is right-veering if and only if the
fractional Dehn twist coefficient is $c\geq 0$. (Observe that $c=0$
corresponds to the identity diffeomorphism.) There is a short list
of periodic elements in $SL(2,\Z)$, up to conjugation:
$$A_1= \begin{pmatrix} 0 & 1 \\ -1 & 0\end{pmatrix}, ~~A_2=
\begin{pmatrix} 1 & 1 \\ -1 &0\end{pmatrix}, ~~A_3 =\begin{pmatrix} 0 &
1 \\ -1 & -1\end{pmatrix},$$ together with $-A_1, -A_2,-A_3$. The
least right-veering lifts (they are right-veering, but leftmost
amongst right-veering lifts) of $A_i$ are given by
$a_1=\sigma_1\sigma_2\sigma_1$, $a_2=\sigma_1\sigma_2$, and
$a_3=(\sigma_1\sigma_2)^2$, respectively. These correspond to
``rotations'' by amounts ${\pi\over 2}, {\pi\over 3}, {2\pi\over 3}$
in the clockwise direction. (Hence $c$ must be in multiples of
${1\over 4}$ or ${1\over 6}$.) The least right-veering lifts of
$-A_i$ are given by multiplying the above lifts $a_i$ by the central
element $(\sigma_1\sigma_2\sigma_1)^2$, and the other right-veering
lifts of $\pm A_i$ are $a_i(\sigma_1\sigma_2\sigma_1)^{2k}$, $k\geq
0$. Notice that all the right-veering lifts are products of positive
Dehn twists. Since any right-veering periodic $h$ is conjugate to
one of the above lifts, $h$ must also be a product of positive Dehn
twists. Hence, if $h$ is periodic, then $h\in \veer$ if and only if
$h\in \dehn$.   In other words, there is no difference between
$\veer$ and $\dehn$ for periodic elements.

Finally, if $h$ is reducible, then $h$ can be written as
$(\sigma_1\sigma_2 \sigma_1)^{2n} R_{\gamma}^{m}$, where $m$ and $n$
are integers, and $R_{\gamma}$ is a positive Dehn twist about some
nonseparating curve $\gamma$. By Corollary~3.4 of \cite{HKM2} and
the fact that $\veer\subset \dehn$ are monoids, we see that $h$ is
right-veering if and only if either $n>0$, or $n=0$ and $m\geq 0$.
Later we will show that if $n>0$ and $m\ll 0$, then $h$ is
right-veering but not a product of positive Dehn twists.

\subsection{The Rademacher function} \label{rade}
\begin{figure}[ht]
\begin{overpic}[height=10cm]{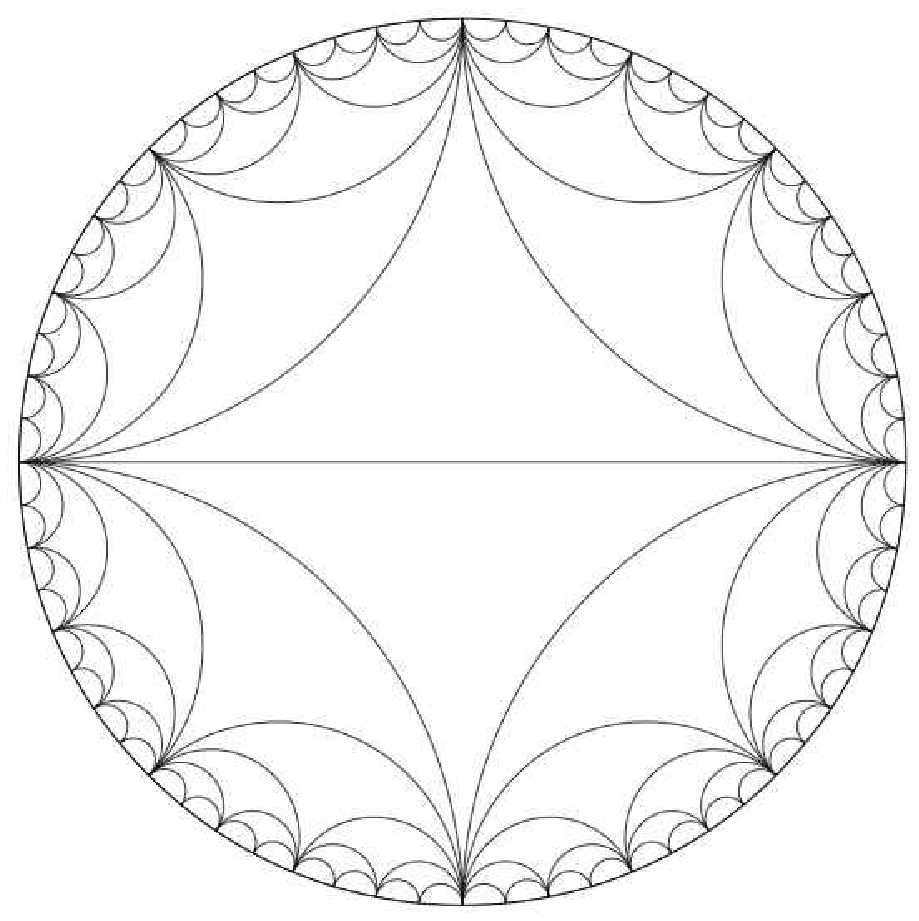}
\put(86,38){$\dfrac{0}{1}$} \put(77,65){$\dfrac{1}{2}$}
\put(48.5,77){$\dfrac{1}{1}$} \put(20,65){$\dfrac{2}{1}$}
\put(11,38){$\dfrac{1}{0}$} \put(18,10){$-\dfrac{2}{1}$}
\put(47,-.8){$-\dfrac{1}{1}$} \put(73.5,10){$-\dfrac{1}{2}$}
\put(50,40){$\scriptstyle 0$} \put(38,51){$\scriptstyle -1$}
\put(39,61){$\scriptstyle 0$} \put(24,55){$\scriptstyle -1$}
\put(26.5,51){$\scriptstyle -2$} \put(41.5,29){$\scriptstyle 1$}
\put(56,29){$\scriptstyle -1$} \put(60,51){$\scriptstyle 1$}
\put(71,51){$\scriptstyle 2$} \put(60,61){$\scriptstyle 0$}
\end{overpic}
\caption{The Farey tessellation and values of the Rademacher
function on the tessellation.} \label{farey-rademacher}
\end{figure}

Consider the action of $PSL(2,\Z)$ on the upper half-plane $\H^2$
and hence on the Farey tessellation of the unit disk $D^2$. Given
$\begin{pmatrix} a & b \\ c& d\end{pmatrix}\in PSL(2,\Z)$, it acts
on $\H^2$ by mapping $z\mapsto {c+dz\over a+bz}$. In particular, if
$z={p\over q}$ is a rational point on the $x$-axis, then ${p\over q}
\mapsto {cq+dp\over aq+bp}$. Under the
correspondence ${p\over q} \leftrightarrow \begin{pmatrix}q \\
p\end{pmatrix}$, the action of $\begin{pmatrix} a & b \\ c&
d\end{pmatrix}$ on ${p\over q}$ is given by left multiplication on
$\begin{pmatrix}q \\ p\end{pmatrix}$. Figure~\ref{farey-rademacher}
shows the Poincar\'e disk model for $\H^2$ with the points on
$S^1_{\infty}$ labeled with the corresponding points on the $x$-axis
of the upper half-plane model.

We now define the {\em Rademacher function} $\Phi:
PSL(2,\Z)\rightarrow \Z$. Much of what follows is taken from
\cite{BG,GG1,GG2}. First observe that $PSL(2,\Z)$ is isomorphic to
the free product $\Z/2\Z * \Z/3\Z$, whose generators are $\pm A$ and
$\pm B$. Hence any element $g\in PSL(2,\Z)$ is uniquely written as
$B^{r_1}AB^{r_2}\dots B^{r_k}$, where $r_1, r_k=0,1,$ or $-1$ and
$r_i=-1$ or $1$, otherwise. We then define $\Phi(g)= \sum_{i=1}^{k}
r_i$.

For a more geometric interpretation of $\Phi$, we will describe how
$\Phi$ can be viewed as a function on the set of edges of the Farey
tessellation. (See Figure~\ref{farey-rademacher}.)  We use the
bijection between directed edges of the Farey tessellation and
$PSL(2,\Z)$, in which a directed edge $a\rightarrow b$ is identified
with the element $g\in PSL(2,\Z)$ which maps the slope $0$ to $a$
and the slope $\infty$ to $b$.  (In other words, $g$ is an
orientation-preserving linear map which sends $(1,0)$ to a shortest
integral vector with slope $a$ and $(0,1)$ to a shortest integral
vector with slope $b$.) Notice that if $g$ corresponds to
$a\rightarrow b$, then $gA$ corresponds to $b\rightarrow a$.  Since
right multiplication by $A$ does not change the value of $\Phi$, it
follows that $\Phi$ induces a function on the set of (undirected)
edges of the Farey tessellation.  Undirected edges will be written
as $ab$.

Again with $g$ corresponding to $a\rightarrow b$, choose $c$ so that
$a$, $b$, and $c$ form the vertices of a triangle in {\em clockwise}
order. Then $gB$ corresponds to $c\rightarrow a$ and $gB^{-1}$
corresponds to $b\rightarrow c$.  Since the value of $\Phi$ on the
identity map (or, equivalently, the edge $0 \rightarrow \infty$) is
$0$, the value on the edge corresponding to any $g$ can be computed
as follows.  Let $p$ be a point on the edge $0\rightarrow \infty$.
Then $\Phi(g)$ equals the number of right turns minus the number of
left turns for a geodesic from $p$ to $g(p)$. Here, a {\em right
turn} refers to an arc from the edge $ab$ to $ac$ and a {\em left
turn} refers to an arc from $ab$ to $bc$.

\s\n {\bf Remark.} In our definition of $\Phi$, we count the number
of right turns minus the number of left turns.  In \cite{BG,GG1},
the authors count the number of left turns minus the number of right
turns.  The definitions of $\Phi: PSL(2,\Z)\rightarrow \Z$ agree (at
least on the set of hyperbolic elements), and the discrepancy is due
to the difference in defining the action on $\H$. (If we defined
$z={z_1\over z_2}$ instead of $z={z_2\over z_1}$, then we would be
counting left turns minus right turns.)

\s
One easily observes that $\Phi: PSL(2,\Z)\rightarrow \Z$ is a {\em
quasi-morphism}. A {\em quasi-morphism} is a map $\phi: G\rightarrow
A$, where $G$ is a group and $A=\Z$ or $\R$, together with a constant
$C$, such that $|\phi(g_1g_2)- \phi(g_1)-\phi(g_2)|\leq C$ for all
$g_1,g_2\in G$. $\Phi$ is also not quite a homomorphism, as can be
seen by taking $g_1$ which ends with $B$ and $g_2$ which begins with
$B$.

\subsection{The linking number and rotation number} \label{link}

Let $B_n$ be the braid group on $n$ strands. Then the {\em linking
number} is a homomorphism $lk: B_n\rightarrow \Z$, defined as
follows: if we write $\sigma\in B_n$ as $\sigma_{i_1}^{j_1}
\dots\sigma_{i_k}^{j_k}$, where $\sigma_1,\dots, \sigma_{n-1}$ are
the standard positive half-twists that generate $B_n$, then
$lk(\sigma) = j_1+\dots+j_k$. The linking number $lk$ is a
homomorphism because $B_n$ has relations only of the type
$\sigma_i\sigma_{i+1}\sigma_i=\sigma_{i+1}\sigma_i \sigma_{i+1}$ and
$\sigma_i\sigma_j=\sigma_j\sigma_i$, i.e., those that leave the sums
of the exponents constant.  In fact, it is the {\em unique}
homomorphism $B_n\rightarrow \Z$ (up to a constant multiple).

There is another invariant of $Aut(S,\bdry S)\simeq B_3$, which we
will call the {\em rotation number} $\rot(h)$, which roughly
measures the number of times $h$ rotates around $\bdry S$.  The
normalization is such that $\rot(R_{\bdry S})=1$.  Just as there is
no homomorphism $Aut(T)\rightarrow Aut(S,\bdry S)$ which splits
Equation~\ref{exact}, there are non-canonical choices involved in
our definition of $\rot(h)$.

Let $h\in Aut(S,\bdry S)$.  We write $\sigma$ for the corresponding
element in $B_3$, and $\overline{\sigma}$ or $a\rightarrow b$ for
its image in $PSL(2,\Z)$.  We consider four cases. (If $a\rightarrow
b$ is $0\rightarrow \infty$ (resp.\ $\infty\rightarrow 0$), then it
is defined in Cases 1 and 4 (resp.\ Cases 2 and 3), and the two
definitions agree.)

\s\n {\bf Case 1.} If $0\leq a<b\leq +\infty$, then we claim that
$\sigma$ can be uniquely written as $(\sigma_1\sigma_2\sigma_1)^{2n}
w$, where $w$ is a word generated by $\sigma_1^{-1}$ and $\sigma_2$
and no inverses of these are allowed.  First observe that
$BA=\overline \sigma_2$ and $B^{-1}A= \overline
\sigma_1^{-1}(\overline\sigma_1
\overline\sigma_2\overline\sigma_1)^2 = \overline \sigma_1^{-1}$.
Since $0\leq a<b\leq +\infty$, $a\rightarrow b$ can be written as
$B^{\pm 1} A B^{\pm 1} A\dots B^{\pm 1} A$. (Recall the
interpretation of $gA$, $gB$, and $gB^{-1}$ from Section~\ref{rade},
where $g\in PSL(2,\Z)$ is viewed as a directed edge. The Farey
tessellation and the dual graph indicate how to move the edge
$0\rightarrow \infty$ to $a\rightarrow b$.) This implies that
$\overline{\sigma}$ can be written as a word $\overline w$ in
$\overline{\sigma}_1^{-1}$ and $\overline{\sigma}_2$.  The element
$\sigma$ can therefore be written as
$(\sigma_1\sigma_2\sigma_1)^{2n}w$, by observing that the kernel of
the projection $B_3\rightarrow PSL(2,\Z)$ is generated by the
central element $(\sigma_1\sigma_2\sigma_1)^2$. The uniqueness is a
consequence of fact that the dual graph to the Farey tessellation is
a tree (and hence there is a unique geodesic between any two
vertices of the graph).

\s\n {\bf Case 2.} If $-\infty\leq a<b\leq 0$, then $\sigma$ can
uniquely be written as $(\sigma_1\sigma_2\sigma_1)^{2n+1}w$, where
$w$ is a word in $\sigma_1^{-1}$ and $\sigma_2$.  Here,
$a\rightarrow b$ can be written as $A(B^{\pm 1} A B^{\pm 1} A\dots
B^{\pm 1} A)$.

\s\n {\bf Case 3.} If $0\leq b<a\leq +\infty$, then $\sigma$ can
uniquely be written as $(\sigma_1\sigma_2\sigma_1)^{2n+1}w$, where
$w$ is a word in $\sigma_1$ and $\sigma_2^{-1}$.  This is because
$a\rightarrow b$ can be written as $B^{\pm 1}(AB^{\pm 1} AB^{\pm 1}
\dots AB^{\pm 1})$ and $AB=\overline{\sigma}_1$, $AB^{-1}=
\overline{\sigma}_2^{-1}$.   Moreover,
$B=(\overline\sigma_1\overline
\sigma_2\overline\sigma_1)^{-1}\overline\sigma_1$ and $B^{-1}
=(\overline\sigma_1\overline\sigma_2\overline\sigma_1)\overline
\sigma_2^{-1} =
(\overline\sigma_2\overline\sigma_1\overline\sigma_2)
\overline\sigma_2^{-1}=\overline\sigma_2\overline\sigma_1$.

\s\n {\bf Case 4.} If $-\infty\leq b<a\leq 0$, then $\sigma$ can
uniquely be written as $(\sigma_1\sigma_2\sigma_1)^{2n}w$, where $w$
is a word in $\sigma_1$ and $\sigma_2^{-1}$.  Indeed, $a\rightarrow
b$ can be written as $AB^{\pm 1} AB^{\pm 1}\dots AB^{\pm 1}$.

\s We now define the {\em rotation number $\rot(\sigma)$} to be
${k\over 4}$, where $\sigma=(\sigma_1\sigma_2\sigma_1)^k w$ as
above.

\begin{thm}\label{lk}
$lk(\sigma) = 12 ~ \rot(\sigma) + lk (w) = 12
~\rot(\sigma)+\Phi(\overline{\sigma}).$
\end{thm}

\begin{proof}
Since $lk$ is a homomorphism and
$\sigma=(\sigma_1\sigma_2\sigma_1)^{4rot(\sigma)}w$ by the
definition of the rotation number, the first equality follows.  To
see that $lk(w)=\Phi(\overline{\sigma})$, first note that
$\overline{\sigma}_1^{-1}=B^{-1}A$ and $\overline{\sigma}_2 = BA$.
If $w$ is a word in $\sigma_1^{-1}$ and $\sigma_2$, as is the case
in Case 1, then the corresponding word in $A$ and $B^{\pm 1}$
involves no cancellation of powers of $B$. It follows that $lk(w)$,
the exponent sum of the $\sigma_1$'s and $\sigma_2$'s, is the same
as the exponent sum of the $B$'s in the word corresponding to $w$;
this, by definition, is $\Phi(\overline{w})$. Finally, since
$\overline{\sigma}_1\overline{\sigma}_2\overline{\sigma}_1=(AB)(BA)(AB)=A$,
powers of
$\overline{\sigma}_1\overline{\sigma}_2\overline{\sigma}_1$
contribute nothing to the Rademacher function and therefore
$\Phi(\overline{w})=\Phi(\overline{\sigma})$.
\end{proof}

\subsection{$\veer$ vs. $\dehn$}
In this subsection we prove Theorem~\ref{torusmain} and
Lemma~\ref{rotationlarge}, which together comprise Theorem~\ref{thm:
veer not in dehn}, and explore some consequences.

Observe that the linking number is useful in detecting braids which
are not quasipositive:

\begin{lemma} \label{notrv}
If one of the following holds, then $\sigma\in B_n$ is not
quasipositive:
\be
\item $lk(\sigma)<0$.
\item $lk(\sigma)=0$ but $\sigma\not=1$. \item $lk(\sigma)=1$ and
$\sigma$ is not conjugate to a half-twist.
\ee
\end{lemma}

\begin{proof}
This follows immediately from the fact that $lk$ is a homomorphism
and consequently is constant on conjugacy classes.
\end{proof}

Theorem~\ref{lk} and (1) of Lemma~\ref{notrv} together imply that if
$-\Phi(\overline\sigma) >12~\rot(\sigma)$, then $\sigma$ is not
quasipositive. We can actually do better:

\begin{thm}\label{torusmain}
If $-\Phi(\overline \sigma) \geq 10~ \rot(\sigma)$, then $\sigma$ is
not quasipositive.
\end{thm}

\begin{proof}
We claim that if $\sigma =\sigma''\sigma'$, where $\sigma''$ is a
positive Dehn twist, then the triple
$(lk(\sigma)-lk(\sigma'),\rot(\sigma)-\rot(\sigma'),
\Phi(\overline\sigma)-\Phi(\overline\sigma'))$ is one of $(1,0,1)$,
$(1,{1\over 4}, -2)$, or $(1,{1\over 2}, -5)$; moreover, if
$\sigma'=id$, then only $(1,0,1)$ and $(1,{1\over 4},-2)$ are
possible. We then observe that $-\Phi(\overline{\sigma})\leq
10~\rot(\sigma)$ if $\sigma$ is a product of positive Dehn twists.
Since the first positive Dehn twist only contributes $(1,0,1)$ or
$(1,{1\over 4},-2)$, we find that the strict inequality
$-\Phi(\overline{\sigma})=10~\rot(\sigma)$ is never attained.

The proof of the claim is a case-by-case analysis.  Suppose
$\overline\sigma$ is written as $a\rightarrow b$ and
$\overline\sigma'$ as $a'\rightarrow b'$. The slope of the Dehn
twisting curve of $\sigma''$ is denoted by $c$.  To visualize the
action of this Dehn twist, consider the infinite collection of arcs
of the Farey tessellation which end at $c$.  The Dehn twist fixes
the point $c$ and maps each arc clockwise to the next arc. Observe
that the three cases below are sufficient, by reversing arrows or by
taking negatives if necessary.

\s\n {\bf Remark.}  The claim is intuitively reasonable if we
consider the ``amount of rotation about $\bdry S$'' effected by each
positive Dehn twist.  The difficulty is that this quantity has a
precise yet non-canonical meaning.  We instead choose to keep track
of $\Phi$, and the method of proof will be useful later in
Section~\ref{examples}.

\s\n {\bf Case 1.}  Suppose $a'\rightarrow b'$ is $0\rightarrow
\infty$.

If $c=\infty$, then $\overline\sigma$ is $-1\rightarrow \infty$, and
if $c=0$, then $\overline\sigma$ is $0\rightarrow 1$. In both cases
$\Phi$ changes by $+1$, and hence $\rot$ by $0$, in view of
Theorem~\ref{lk}.

If $0<c<+\infty$, then let $c,d,e$ be the vertices of a triangle of
the Farey tessellation in clockwise order, so that $a'b'$ and $de$
are in the same connected component of $D^2$ cut open along $cd$ and
$ce$. It could happen that $a'b'=de$. If we apply a positive Dehn
twist about $c$, then $cde$ will be mapped to the adjacent triangle
$cef$. Let $P$ be a word in $L$ and $R$ which records the left and
right turns taken on the geodesic from $0\infty$ to $de$.  (For
example, $LRLL$ means you first take a left turn and then a right
turn, followed by two left turns.) If $P^{-1}$ is obtained from $P$
by reversing the word order and changing an $R$ to an $L$ and an $L$
to an $R$ (for example, if $P=LRLL$, then $P^{-1}=RRLR$), then the
path from $0 \infty$ to $ab$ is given by $PLLP^{-1}$. See
Figure~\ref{farey-onedehn}. $\Phi$ changes by $-2$ and hence $\rot$
by ${1\over 4}$. Similarly, if $-\infty<c<0$, then $\Phi$ changes by
$-2$ and $\rot$ by ${1\over 4}$.

\begin{figure}[ht]
\s
\begin{overpic}[height=7cm]{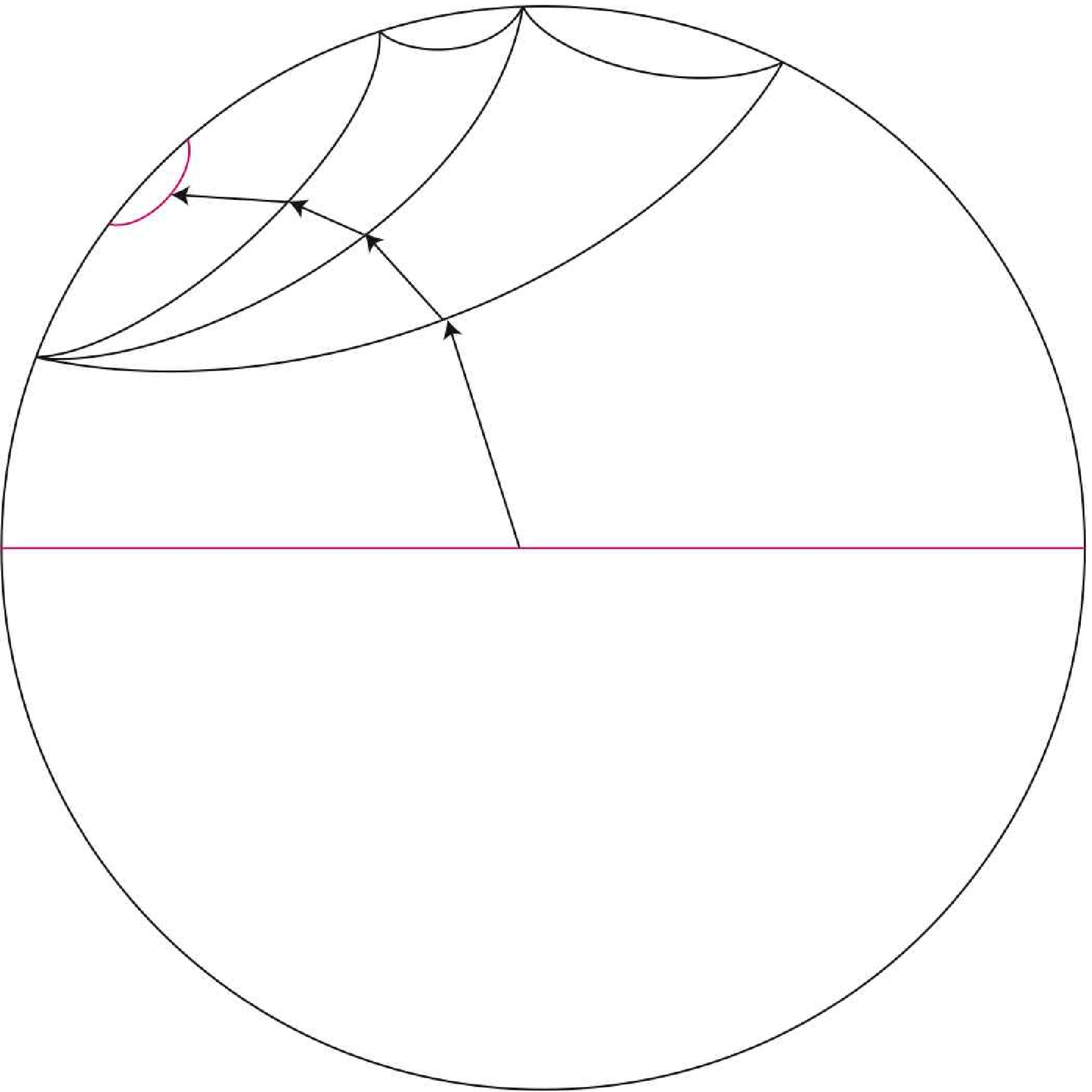}
\put(-19,48.6){\small{$b'=\infty$}}
\put(101.5,48.67){\small{$a'=0$}}
\put(47,101.7){\small{$c$}}
\put(71.8,96.8){\small{$d$}}
\put(-1.45,67.4){\small{$e$}}
\put(33,100.4){\small{$f$}}
\put(15,90){\small{$b$}}
\put(4.25,78.7){\small{$a$}}
\put(46,60){\tiny{$P$}}
\put(38.4,75.4){\tiny{$L$}}
\put(31,82){\tiny{$L$}}
\put(19.5,83.5){\tiny{$P^{-1}$}}
\end{overpic}
\caption{The path $PLLP^{-1}$ from $0\infty$ to $ab$ in Case 1.}
\label{farey-onedehn}
\end{figure}

\s\n {\bf Case 2.}  Suppose $0\leq a'< b'<+\infty$.

If $c=a'$, then $\Phi$ changes by $+1$.  If $c=b'$, then $\Phi$
changes by $+1$ if $0\leq a<a'$, and by $-2$ if $b'< a\leq +\infty$.
If $a'<c<b'$, then $\Phi$ changes by $-2$ as in Case 1.

If $0\leq c< a'$, then let $cde$ be as above, i.e., such that $a'b'$
and $de$ are in the same connected component of $D^2$ cut open along
$cd$ and $ce$,  and let $cef$ be the image of $cde$ under the
positive Dehn twist about $c$. Suppose first that $d\geq b'$. If we
draw a diagram like the one in Figure~\ref{farey-onedehn}, we can
see that the path from $0\infty$ to $a'b'$ can be labeled by $P'LQ$,
where $P'$ is the path from $0\infty$ to $cd$, $L$ is the left turn
at $d$, i.e., from $cd$ to $de$, and $Q$ is the path from $de$ to
$a'b'$.  Then the path from $0\infty$ to $ab$ is given by $P'RLQ$.
Here $R$ is a right turn around $c$ taking $dc$ to $ce$ and $L$ is a
right turn around $e$ taking  $ce$ to $ef$. We easily see that in
this case $\Phi$ changes by $+1$. Next suppose that $0\leq d<c$.  If
$P_1$ is the path from $0\infty$ to $a'b'$, and $P_2$ is the path
from $a'b'$ to $de$, then we can write the paths as $P_1=P_1'LQ$ and
$P_2=Q^{-1}L P_2'$. Then the path from $0\infty$ to $ab$ is given by
$\widehat{P_1P_2}LL P_2^{-1}$, where $\widehat{\mbox{ }\mbox{ }}$
indicates a contraction. More precisely, since $P_1=P_1'LQ$ and
$P_2=Q^{-1}L P_2'$, where the two $L$'s are around two vertices of
the same triangle,  we get $\widehat{P_1P_2}=P_1'RP_2'$, where $R$
corresponds to the right turn at the third vertex. In the end the
path from $0\infty$ to $ab$ is  $P_1'RP_2'LL P_2^{-1}=
P_1'RP_2'LL(P'_2)^{-1}RQ$. Recalling that $P_1=P_1'LQ$,  we see that
$\Phi$ changes by $+1$ as well.

Suppose $-\infty < c<0$. If the path from $0\infty$ to $a'b'$ is
$P_1$, then the path from $0\infty$ to $ab$ is of the form $P_2 LL
P_2^{-1} P_1$.  Therefore, $\Phi$ changes by $-2$.

Next suppose $b'<c < +\infty$.  If $a'\geq e > f\geq 0$, then the
path from $0\infty$ to $a'b'$ is given by $P_1LRP_2$, where $L$ and
$R$ are turns through $cef$ and $cde$. Then the path from $0\infty$
to $ab$ is given by $P_1RP_2$, and $\Phi$ changes by $+1$. If
$a'\geq e\geq 0$ and $f>c$, then $a,b$ satisfies one of the
following: (i) $e\geq a>b \geq 0$, (ii) $-\infty \leq a,b \leq 0$,
(iii) $a,b \geq f$, or (iv) $e\geq a\geq 0$ and $b>f$.  In any case,
we can write $P_1RRP_2$ for the path from $0\infty$ to $a'b'$, where
$P_1$ is the path from $0\infty$ to $ef$, the $R$'s rotate about $e$
and $P_2$ is the path from $de$ to $a'b'$. In case  (i), the word
$P_1RRP_2$ for the path from $0\infty$ to $a'b'$ is transformed to
$\widehat{P_1P_2}$. More precisely, we can write $P_1=P_1'LQ$ and
$P_2=Q^{-1}L P_2'$, and then we get $P_1'R P_2'$. Thus $\Phi$
changes by $+1$.  In (ii), the word $P_1RRP_2$ is transformed to
$P_2'$, where $P_2=P_1^{-1}P_2'$. $\Phi$ changes by $-2$.  In (iii),
$P_1RRP_2$ is transformed to $\widehat{P_1P_2}$, where $P_1=P_1'RQ$
and $P_2=Q^{-1}RP_2'$ and $\widehat{P_1P_2}=P_1'LP_2'$. This time
$\Phi$ changes by $-5$. In (iv), $P_1P_2RRP_2^{-1}$ is transformed
to $P_1$, and $\Phi$ changes by $-2$. If $e > f>c$, $P_1$ is the
path from $0\infty$ to $a'b'$ and $P_2$ is the path from $a'b'$ to
$de$, then the path from $0\infty$ to $ab$ is
$\widehat{P_1P_2}LLP_2^{-1}$, and $\Phi$ changes by $-5$.

Finally take $c= \infty$.  If $cde=\infty 1 0$, then $RP$ maps to
$LP$ and $\Phi$ changes by $-2$.   Otherwise, $P_1LRP_2$ maps to
$P_1RP_2$ and $\Phi$ changes by $+1$.

\s\n {\bf Case 3.} Suppose $a'\rightarrow b'$ is $a'\rightarrow
\infty$, where $a'$ is a nonnegative integer.

If $c=\infty$ or $c=a'$, then $\Phi$ changes by $+1$.   If
$a'<c<+\infty$, then $\Phi$ changes by $-2$ as in Case 1. If $0\leq
c<a'$, then $\Phi$ changes by $+1$ as in Case 2.  Finally, if
$-\infty< c<0$, then $\Phi$ changes by $-2$ as in Case 2.  Notice
that in this case $\Phi$ cannot change by $-5$.
\end{proof}

Theorem~\ref{torusmain} is effective when used in conjunction with
the following lemma:

\begin{lemma}\label{rotationlarge}
If $rot(\sigma)\geq {1\over 2}$, then $\sigma$ is right-veering.
\end{lemma}

Before proceeding with the proof, we briefly discuss the action of
$h\in Aut(S,\bdry S)$ on the universal cover $\widetilde S$ of $S$,
as described in \cite{HKM2}. In this paragraph we assume that the
Euler characteristic $\chi(S)$ is negative, i.e., $S$ is not a disk
or an annulus.  Endow $S$ with a hyperbolic metric for which $\bdry
S$ is geodesic. The universal cover $\pi: \tilde S\rightarrow S$ can
then be viewed as a subset of the Poincar\'e disk $D^2 =
\mathbb{H}^2 \cup S^1_\infty$. Now let $L$ be a component of
$\pi^{-1}(\bdry S)$. If $h\in Aut(S,\bdry S)$, choose a lift $\tilde
h$ of $h$ that is the identity on $L$. The closure $\tilde S$ in
$D^2$ is geodesically convex.  One portion of $\bdry \tilde S$ is
$L$ and the complement of the closure of $L$ in $\bdry \tilde S$
will be denoted $L_\infty$. Note that $L_\infty$ is homeomorphic to
$\R$. Orient $L_\infty$ using the boundary orientation of $\tilde S$
and then linearly order $L_\infty$ so that moving in an
orientation-preserving sense increases the order. The lift $\tilde
h$ induces a homeomorphism $h_\infty:L_\infty \to L_\infty$. By
Theorem~2.2 of \cite{HKM2}, $h$ is right-veering if and only if $z
\ge h_\infty(z)$ for all $z\in L_\infty$.

\begin{proof}[Proof of Lemma~\ref{rotationlarge}.]
This can be proved on a case-by-case basis, as in the definition of
the rotation number in Section~\ref{link}.  We will treat Case 1,
i.e., $0\leq a<b\leq +\infty$, and leave the other cases to the
reader.  Since $rot(\sigma)\geq {1\over 2}$, we have
$\sigma=(\sigma_1\sigma_2\sigma_1)^{2n}w$, where $n\geq 1$. It
suffices to verify the lemma for $n=1$. In the rest of the proof we
write $\sigma$ if we mean an element in $B_3$, and write $h$ to
denote the corresponding element in $Aut(S,\bdry S)$.

Using the notation from the paragraph preceding the proof, pick a
basepoint $x\in \bdry S$ and a lift $\widetilde{x}\in L$.  We can
endow $L_\infty$ with a nondecreasing continuous function $\theta:
L_\infty \rightarrow \R$ so that any properly embedded, oriented arc
$\alpha:[0,1]\rightarrow S$ with $\alpha(0)=x$ and slope $s$
satisfies $\theta(\widetilde\alpha(1))\equiv \theta_s(\mbox{mod }
2\pi)$, where $\widetilde\alpha$ is a lift of $\alpha$ to
$\widetilde{S}$ whose initial point is $\widetilde x$ and $\theta_s$
is the standard angle that a line of slope $s$ makes with a line of
slope $0$. (Here $\theta$ would be an angular coordinate on
$K_\infty$, obtained from $L_\infty$ by quotienting each connected
component of $\pi^{-1}(\bdry S)$ besides $L$ to a point.)

Let $\alpha$ and $\beta$ be properly embedded, oriented arcs based
at $x$ with slopes $a$ and $b$, such that $0\leq
\theta(\widetilde{\alpha}(1))< \theta(\widetilde{\beta}(1))\leq
{\pi\over 2}$. The element $h_\infty$ maps the interval
$[0,{\pi\over 2}]$ to
$[\theta(\widetilde{\alpha}(1))-\pi,\theta(\widetilde{\beta}(1))-\pi]$
and maps the interval $[-{\pi\over 2},0]$ to
$[\theta(\widetilde{\beta}(1))-2\pi,\theta(\widetilde{\alpha}(1))-\pi]$.
By applying the same argument to other intervals, we see that
$\sigma$ is right-veering. The other cases are similar.
\end{proof}

\n {\bf Remark.}  Observe that, in order to show that $\sigma$ is
right-veering, it is not sufficient to verify that two properly
embedded arcs of $S$ corresponding to an integer basis of $\Z^2$ get
mapped to the right.

\s\s In the rest of the subsection we give some consequences of the
above discussion.

\begin{cor}   \label{infmany}
For the punctured torus $S$, there are infinitely many pseudo-Anosov
diffeomorphisms $h\in \veer$ with arbitrarily large fractional Dehn
twist coefficients $c$, which are not in $\dehn$.
\end{cor}

\begin{proof}
As before, we switch freely between $\sigma\in B_3$ and its
corresponding $h\in Aut(S,\bdry S)$. Choose $\sigma=(\sigma_1
\sigma_2 \sigma_1)^{2n}w$ with $n\ge 1$ as in Case~1, but with
$0<a<b<+\infty$.  Then, the action of $\overline\sigma\in PSL(2,\Z)$
on the circle at infinity $S^1_{\infty}$ of the Farey tessellation
has two fixed points. Therefore $\overline\sigma$ is Anosov and
$h\in Aut(S,\bdry S)$ is pseudo-Anosov.  (Alternatively, one can
compute the trace of $\overline{\sigma}$, and show that it is $>2$
or $<-2$, since the entries are all positive or all negative.) Since
$n\ge 1$, all such $h$ are right-veering by
Lemma~\ref{rotationlarge}. On the other hand, if $w$ is chosen so
that $\#(\mbox{$\sigma_1^{-1}$ terms})- \#(\mbox{$\sigma_2$ terms})$
is sufficiently large, then $h\notin Dehn^+(S, \bdry S)$.
\end{proof}

To rephrase Corollary~\ref{infmany} in terms of the braid group
$B_n$, we recall Thurston's {\em left orderings} of $B_n$. Let $S$
be the double branched cover of the unit disk, branched along $n$
points. Thurston defined {\em left orderings} of $Aut(S, \partial
S)$ (and hence the left orderings on $B_n$) via the double branched
cover $S$: Fix $z\in L_\infty$. Given $h,g\in Aut(S,
\partial S)$, define $h\ge_z g$ if $h(z) \ge g(z)$. Such an
ordering is called a {\em left ordering} because it preserved by
left multiplication. (Of course, $\ge_z$ may not be a {\em total
order}, but that is not an important issue here.) The following is a
rephrasing of Corollary~\ref{infmany}.

\begin{cor}
There exist infinitely many pseudo-Anosov braids $\sigma\in B_3$ for
which $id \geq \sigma$ using any of the left orderings of $B_3$
defined by Thurston, but which are not quasipositive.
\end{cor}

\n {\bf Example:} Let
$\sigma=(\sigma_1\sigma_2\sigma_1)^2\sigma_1^{-m}$. Then
$\overline\sigma = \begin{pmatrix} -1& 0 \\ -m& -1\end{pmatrix}$.
This is the right-veering lift of $\overline\sigma\in SL(2,\Z)$ with
the least amount of rotation to the right. By
Theorem~\ref{torusmain}, if $m\geq 5$, then $\sigma$ is not
quasipositive. On the other hand, we claim that, for $m\leq 4$,
$\sigma$ can be written as a product of positive Dehn twists.  It
suffices to prove the claim for
$(\sigma_1\sigma_2\sigma_1)^2\sigma_1^{-4}$.  Indeed using the fact
that:
$$(\sigma_1\sigma_2\sigma_1)\sigma_1=\sigma_2(\sigma_1\sigma_2\sigma_1),
~~(\sigma_1\sigma_2\sigma_1)\sigma_2=\sigma_1(\sigma_1\sigma_2\sigma_1),$$
we write:
$$(\sigma_1\sigma_2\sigma_1)^2\sigma_1^{-4}
=(\sigma_1\sigma_2\sigma_1)(\sigma_1\sigma_2\sigma_1^{-1})\sigma_1^{-2}
=(\sigma_2\sigma_1\sigma_2^{-1})(\sigma_1\sigma_2\sigma_1)\sigma_1^{-2}
=(\sigma_2\sigma_1\sigma_2^{-1})(\sigma_1\sigma_2\sigma_1^{-1}).$$

More generally, we can show that
$(\sigma_1\sigma_2\sigma_1)^{2n}\sigma_1^{-m}$ is quasipositive if
$m\leq 4n$ but not quasipositive if $m\geq 5n$. In general, we do
not know what happens for $m$ strictly between $4n$ and
$5n$.\footnote{The referee has informed us that the following holds:
For any $m$, the least $\Psi(m)$ for which
$(\sigma_1\sigma_2\sigma_1)^{\Psi(m)}\sigma_1^{-m}$ is quasipositive
is $\Psi(m)=2k+1$ for $m=5k,5k+1,5k+2$, and $\Psi(m)=2k+2$ for
$m=5k+3,5k+4$.  In each case,
$(\sigma_1\sigma_2\sigma_1)^{\Psi(m)-1}$ is not quasipositive by
Theorem~\ref{torusmain}, so Theorem~\ref{torusmain} gives a tight
bound for parabolic elements.}

\begin{cor}
There does not exist a finite set of generators for $Veer(S,\bdry
S)$ over $Aut(S,\bdry S)$, that is, there is no finite collection
$C$ of elements of $Veer(S,\bdry S)$ such that every element of
$Veer(S,\bdry S)$ can be expressed as a product of positive powers
of elements of $C \cup Dehn^+(S,\bdry S)$.
\end{cor}

\begin{proof}
Consider $\sigma(m)=(\sigma_1\sigma_2\sigma_1)^2\sigma_1^{-m}$ with
$m\ge 5$. By the above example, $\sigma(m)\in \veer-\dehn$. The
homeomorphism $h_\infty(m):L_\infty\rightarrow L_\infty$
corresponding to $\sigma(m)$ sends $[0, \pi/2]$ to
$[\tan^{-1}(m)-\pi, \pi/2 - \pi]$. Notice that every angle is
decreased by at most $\pi$, and the only angles that are decreased
by $\pi$ are $\pi/2+k\pi$.

We claim that if $h_1, h_2\in Veer(S,\bdry S)$ and $\sigma(m)=h_1
h_2$, then one of the $h_i$ is $\sigma(m')$ with $m'\geq m$ and the
other is $\sigma_1^{m'-m}$; the corollary then follows immediately
from the claim.  Effectively we are showing that the $\sigma(m)$ are
the least right-veering among right-veering diffeomorphisms which
are not in $\dehn$.

Arguing by contradiction, let $\sigma(m)=h_1h_2$ be such a
factorization. Since $\sigma(m)\notin Dehn^+(S, \bdry S)$, it is not
possible that both $h_1, h_2 \in Dehn^+(S, \bdry S)$. First we claim
that $h_i$ cannot be freely homotopic to a pseudo-Anosov
homeomorphism. Indeed, for a pseudo-Anosov $h_i$ to be
right-veering, it must have fractional Dehn twist coefficient $c\geq
1/2$ by Proposition~3.1 of \cite{HKM2}. Then there is a properly
embedded arc $\alpha$ on $S$ so that
$\theta(\widetilde{\alpha}(1))-\theta(\widetilde{h(m)}(1))>\pi$.
(Take $\alpha$ so that its slope is close to, but slightly larger
than, the stable slope.) Since homeomorphisms which are freely
homotopic to periodic homeomorphisms in $Veer(S,\bdry S)$ are
necessarily in $Dehn^+(S, \bdry S)$, it follows that one of the
factors $h_i$ must be reducible and not in $Dehn^+(S, \bdry S)$.
This means that $h_i$ can be expressed as
$(\sigma_1\sigma_2\sigma_1)^{2n_1} R_\gamma^{n_2}$. Since $h_i$ is
right-veering, $n_1 \ge 0$, but $n_1 = 0$ would imply $n_2 >0$ and
then $h_i\in  Dehn^+(S, \bdry S)$. Also $n_1$ cannot be greater than
or equal to $2$, since the angle of rotation would be too large, and
$h_i$ could not be a factor of $\sigma(m)$. This leaves the
possibility $h_i=(\sigma_1\sigma_2\sigma_1)^{2} R_\gamma^{n_2}$. In
this case the only angles that are decreased by $\pi$ when acted on
by $(\sigma_1\sigma_2\sigma_1)^{2} R_\gamma^{n_2}$ are the angles
corresponding to $\pm \gamma$.  It follows that $\gamma$ has slope
$\infty$; thus $h_i=\sigma(-n_2)$. Letting $m'=-n_2$ and using the
fact that $h_i$ decreases angles by no more than $\sigma(m)$ implies
$m'\geq m$.
\end{proof}

\subsection{An example}   \label{examples}
In this subsection, we will give a computation of an element $h\in
\veer-\dehn$ which does not satisfy the condition of
Theorem~\ref{torusmain}. It is likely that the types of computations
done in the example are amenable to computer calculation, i.e., the
algorithm can probably be done in finite time for the torus.

\s\n {\bf Example:}  $\sigma=(\sigma_1\sigma_2\sigma_1)^2
\sigma_1^{-4}\sigma_2 \sigma_1^{-1}\sigma_2\sigma_1^{-1}$ is in
$\veer-\dehn$.  However, $\Phi=-4$, $\rot={1\over 2}$ and $lk=2$,
and the conditions of Theorem~\ref{torusmain} are not satisfied. Our
strategy is to exploit the fact that $lk=2$, so that $\sigma$ must
be expressed as a product of two positive Dehn twists if
$\sigma\in\dehn$. There are two possibilities: (1) the first Dehn
twist contributes $(1,0)$ to $(\Phi,\rot)$ and the second
contributes $(-5,{1\over 2})$, or (2) the first contributes
$(-2,{1\over 4})$ and the second $(-2,{1\over 4})$.

(1) Referring to the proof of Theorem~\ref{torusmain}, Case 1, the
first Dehn twist sends $0\rightarrow\infty$ to $0\rightarrow 1$ or
$-1\rightarrow\infty$. By looking at the $\sigma$ we are
considering, we see that the image is in the upper half disk of the
Farey tessellation; let us denote the corresponding word by $W$. By
inspecting again the proof of Theorem~\ref{torusmain} we see that if
the second Dehn twist contributes $(-5,{1\over 2})$, then it leaves
$0\rightarrow 1$ or $-1\rightarrow \infty$ in the same half disk
(upper or lower) of the Farey tessellation, so $-1\rightarrow
\infty$ is not possible. In the only possible combination of twists
we easily see that  if $\overline\sigma$ is written as
$a'\rightarrow b'$, then the path $W$ from $0\infty$ to $a'b'$ can
be written as $LPLLP^{-1}L$, where $P$ is some word in $L$ and $R$.
This is a contradiction.

\begin{figure}[ht]
\s
\begin{overpic}[height=5cm]{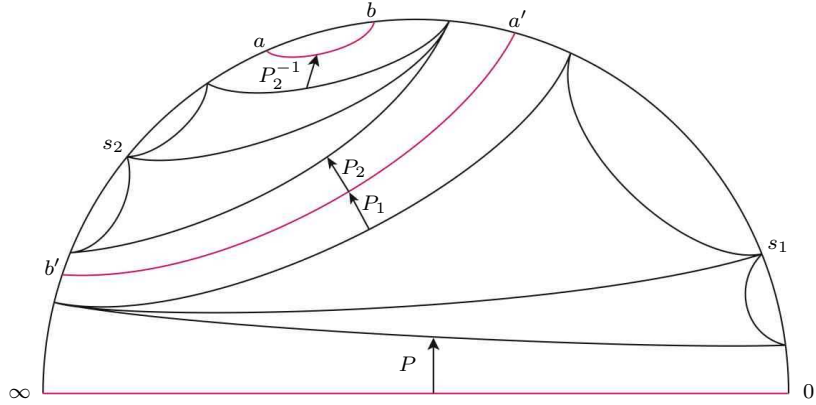}
\put(-4.5,-.5){\tiny{$\infty$}} \put(101.8,-.5){\tiny $0$}
\put(28.3,46.5){\tiny $a$} \put(43.4,50.5) {\tiny $b$}
\put(62.4,49.2){\tiny $a'$} \put(0.3,16){\tiny $b'$}
\put(42.7,24.7){\tiny $P_1$} \put(40,29.5){\tiny $P_2$}
\put(28.9,41.8){\tiny $P_2^{-1}$} \put(47.7,3.2){\tiny $P$}
\put(8.1,32.8){\tiny $s_2$}\put(97,19.5){\tiny $s_1$}
\end{overpic}
\caption{Case (ii)(a) with $s_2>a$. Here $P=P_2LLP_2^{-1}P_1^{-1}$, and
$W= PLLP_1$.} \label{farey-ex}
\end{figure}

(2) The first Dehn twist maps $0\rightarrow \infty$ to $a\rightarrow
b$, where (a) $0< b<a< +\infty$ or (b) $-\infty<a<b<0$. Suppose the
second Dehn twist maps $a\rightarrow b$ to $a'\rightarrow b'$. In
case (a), there are three relevant subcases: (i) the slope $s_2$ of
the second Dehn twist satisfies $b<s_2<a$; (ii) $s_2>a$ and $a'\leq
b$, $b'\geq a$; (iii) $s_2=a$ and $b'>a$. This again follows from
the analysis of the proof of Theorem~\ref{torusmain}. In subcase
(i), $W$ can be written as $P_1LLP_1^{-1}P_2LLP_2^{-1}$. In subcase
(ii), after analyzing all the possible diagrams, we see that the
only relative position of the two adjacent triangles in the Farey
tessellation with vertex $s_2$ that results in a contribution of
$(-2,{1\over 4})$ is the one presented in Figure~\ref{farey-ex}.
Moreover, the edge $a'\rightarrow b'$ can be any edge which
intersects the geodesic from $0\infty$ to the lower triangle with
vertex $s_2$ given in Figure~\ref{farey-ex}.   One particular
possibility for $a'\rightarrow b'$ is given in
Figure~\ref{farey-ex}; this gives $W=P_2 LLP_2^{-1}P_1^{-1}LLP_1$.
The other possibilities for $a'\rightarrow b'$ are edges of the two
triangles with vertex $s_1$ and edges between $0\infty$ and the
lower triangle with vertex $s_1$; they give equations
$W=PLLP^{-1}LL$, $W=LPLLP^{-1}L$, $W=LLPLLP^{-1}$, and
$WP_1=P_1LLP_2LLP_2^{-1}$. In subcase (iii), $W$ can be written as
$LPLLP^{-1}L$. In case (b), there are also two subcases: $s_2>0$ or
$b< s_2\leq 0$. In the former subcase, $W$ can be written as
$P_2LLP_2^{-1}P_1LLP_1^{-1}$. In the latter, we can write $P_1W=
P_2LLP_2^{-1}LLP_1$.

In summary, if $W$ is the path from $0\infty$ to $a'b'$, we need to
show that each of the following equations has no solution:
\begin{eqnarray}
LPLLP^{-1}L &=& W\\
P_1LLP_1^{-1}P_2LLP_2^{-1} &=& W\\
P_2LLP_2^{-1}LLP_1&=&P_1W\\
P_1LLP_2LLP_2^{-1}&=&WP_1
\end{eqnarray}

So far we have only used the facts that $\Phi=-4$, $rot={1\over 2}$
and $lk=2$.  We now show that our specific choice
$W=LLLLRLRL=L^4RLRL$ has no solution to any of the above equations.
The first two equations are immediate. To see that the third
equation has no solution note that we can write it in the form
$XP=PW$ where $P=P_1$. Since $P$ must have the same last letter as
$W$, we can write $P=QY$ and $W=ZY$ (for example, we can take
$Y=RLRL$). Then $XQ=QYZ$. Notice that this is an equation of the
same form as $XP=PW$, but with $Q$ repeated instead of $P$, and that
$Q$ is shorter than $P$ and $YZ$ is a cyclic permutation of $W$.
Notice also that $Q$ is not the empty word since no cyclic
permutation of $W$ is equal to $X$, regardless of the choice of
$P_2$.  This means the argument can be repeated, i.e., $Q$ must have
the same last letter as $Z$ and can be written as $Q=Q_1Y_1$ with
$YZ=Z_1Y_1$, etc. This procedure inductively shortens $Q_i$. Since
we can never reduce to the empty word, this gives us a
contradiction. The fourth equation can be treated in the same way as
the third.

In general, this line of argument can be done for many words. Long
sequences with $RLRLRL...$ are effective, since the equations above
all contain $LL$.

\subsection{Questions} We close this section with some questions.

Suppose $S$ is the punctured torus.  We were able to identify large
swathes of $\veer-\dehn$.  However, we are far from determining all
of $\veer-\dehn$.

\begin{q}
Determine a complete set of invariants that will distinguish
elements of $\veer$ that are not in $\dehn$.
\end{q}

Our initial motivation for undertaking the study of the difference
between $\veer$ and $\dehn$ was to understand the difference between
tight contact structures and Stein fillable contact structures.  The
following question is still open.

\begin{q}
If $h\in \veer-\dehn$ for $S$, then is $(S,h)$ not Stein fillable?
Is it true for a general bordered surface $S$?  Is it true if $S$ is
a punctured torus?
\end{q}

There is some evidence that the answer is yes, which the authors
learned from Giroux. Recall that $(S,h)$ is Stein fillable if and
only if there is some $(S',h')\in \dehn$ so that the open books
$(S,h)$ and $(S',h')$ become the same after performing a sequence of
stabilizations to each (no destabilizations allowed).  The work of
Orevkov~\cite{Or} shows that, in the braid group (or, equivalently,
in the hyperelliptic mapping class group), $\sigma\in B_n$ is
quasipositive if and only if its stabilization in $B_{n+1}$ is
quasipositive.  It is not clear to the authors how to adapt
Orevkov's argument to the more general situation.

By Orevkov, all the $h\in \veer-\dehn$ constructed above for the
punctured torus $S$ are still not quasipositive when stabilized and
viewed in the braid group/hyperelliptic mapping class group.
Moreover, after a certain number of stabilizations, the linking
number of the braid is no longer negative!

\section{Generalizations to the braid group}\label{braidgroup}

We now discuss generalizations of the results from the previous
section to the braid group $B_n$.

One method is to start with $\sigma\in B_3$ which is right-veering
but not in $\dehn$ because $lk<0$, and then embed
$\rho:B_3\hookrightarrow B_n$ (somewhat) canonically by adding extra
strands.  Since $\veer$ is a monoid, one can take products of such
$\rho(\sigma)$, their conjugates in $B_n$, and quasipositive
elements $\sigma'\in B_n$.  Provided the linking number is still
negative, the product is right-veering but not quasipositive.

Another method (presumably slightly more general) is to rephrase the
$lk<0$ condition in terms of the {\em signature} of the braid
closure and the {\em Maslov index}.  This uses, in an essential way,
the work of Gambaudo-Ghys~\cite{GG1,GG2}.  After some preparatory
remarks in Sections~\ref{bounded} and ~\ref{maslov}, we prove
Theorem~\ref{sign+maslov} in Section~\ref{thesignature}.

\subsection{Bounded cohomology}\label{bounded}
In this subsection we interpret Theorem~\ref{torusmain} in terms of
{\em bounded cohomology}.

Since much of the material is probably unfamiliar to specialists in
contact and symplectic geometry, we include a brief summary of
bounded cohomology and the (bounded) Euler class. An excellent
source is \cite{Gh}.

Let $G$ be a group and $A=\Z$ or $\R$.  Then the (ordinary)
cohomology group $H^k(G;A)$ is the cohomology of the chain complex
$(C^k(G;A), \delta)$, where $C^k(G;A)$ is the set of maps $c:
G^{k+1}\rightarrow A$ which are {\em homogeneous}, i.e.,
$c(gg_0,gg_1,\dots, gg_k)= c(g_0,\dots, g_k)$, and $\delta:
C^{k-1}(G;A)\rightarrow C^{k}(G;A)$ is the cochain map:
$$\delta c (g_0,g_1,\dots, g_k)=\sum_{i=0}^{k} (-1)^i c(g_0,\dots, \widehat g_i,
\dots,g_k).$$ The {\em bounded cohomology} group $H^k_b(G;A)$ is the
cohomology of the chain complex $C^k_b(G;A)\subset C^k(G;A)$ of maps
$c: G^{k+1}\rightarrow A$ for which $|\sup_{(g_0,\dots,g_{k})\in
G^{k+1}} c(g_0,\dots,g_{k})|<\infty$. There is a natural map
$H^k_b(G;A)\rightarrow H^k(G;A)$ which is not necessarily injective
or surjective.

The homogeneous cochain  $c: G^{k+1}\rightarrow A$ corresponding to
the inhomogeneous cochain $\overline c: G^k\rightarrow A$ is given
by $c(g_0,\dots,g_k)=\overline c(g_0^{-1}g_1,g_1^{-1}g_2,\dots,
g_{k-1}^{-1}g_k)$. In the other direction, we can find the
inhomogeneous cochain $\overline c$ whose homogenization is $c$ by
setting $\overline c(h_1,\dots,h_k)= c(e,h_1,h_1h_2,h_1h_2h_3,
\dots, h_1h_2 \dots h_k)$. Note that with this dehomogenization, the
coboundary of $C^1(G;A)$ is defined on inhomogeneous maps by
$$\delta \overline c (h_1,h_2)= \overline c(h_1) +\overline c(h_2)
-\overline c(h_1h_2)$$

What we are interested in is $H^2(G;A)$, which classifies isomorphism
classes of central extensions of $G$ by $A$:
\begin{equation}\label{central}
0\rightarrow A\rightarrow \widetilde G\rightarrow G\rightarrow 1.
\end{equation}
The class in $H^2(G;A)$ corresponding to the central extension given
by Equation~\ref{central} is called the {\em Euler class of the
central extension}.

Now define a {\em quasi-morphism} to be a map $\phi: G\rightarrow
A$, together with a constant $C$, such that
$|\phi(g_1g_2)-\phi(g_1)-\phi(g_2)|\leq C$ for all $g_1,g_2\in G$.
Denote by $QM(G;A)$ the $A$-module of quasi-morphisms from $G$ to
$A$. A {\em trivial} quasi-morphism is a quasi-morphism $\phi$ which
is a bounded distance from a genuine homomorphism $\psi$, i.e.,
$\phi-\psi$ is bounded on $G$. (Hence, two quasi-morphisms are
deemed equivalent if their difference is within bounded distance of
a genuine homomorphism.) The following fact can be verified directly
from the definitions.

\s\n
{\bf Fact:} The kernel of $H^2_b(G;A)\rightarrow H^2(G;A)$ is the
quotient of $QM(G;A)$ by the trivial quasi-morphisms.

\subsection{Interpretation of Theorem~\ref{lk} from the viewpoint of bounded
cohomology.}

Let $\mbox{Homeo}_+(S^1)$ be the group of orientation-preserving
homeomorphisms of $S^1$ and $\widetilde{\mbox{Homeo}}_+(S^1)$ be the
universal cover of $\mbox{Homeo}_+(S^1)$.  If we identify
$S^1=\R/\Z$, then an element $\widetilde\gamma$ of
$\widetilde{\mbox{Homeo}}_+(S^1)$ is a periodic
orientation-preserving homeomorphism of $\R$ with period $1$. Define
the {\em translation number} $\Psi: \widetilde{\mbox{Homeo}}_+
(S^1)\rightarrow \R$, where $\Psi(\widetilde\gamma)=
2\widetilde\gamma(0)$ if $\widetilde\gamma(0)\in \Z$ and
$2\lfloor\widetilde\gamma(0)\rfloor+1$ if
$\widetilde\gamma(0)\not\in \Z$. Here $\lfloor\cdot\rfloor$ is the
greatest integer function.  The translation number, roughly
speaking, keeps track of twice the number of times a point is sent
around $S^1$.  The translation number $\Psi$ is a quasi-morphism of
$\widetilde{\mbox{Homeo}}_+(S^1)$, whose coboundary
$$\delta \Psi (g_1,g_2)= \Psi(g_1)+\Psi(g_2)-\Psi(g_1g_2)$$
descends to $\mbox{Homeo}_+(S^1)$ and represents an element in the
second bounded cohomology group $H^2_b(\mbox{Homeo}_+(S^1);\Z)$.

Via the standard action of $PSL(2,\R)$ on $\R\P^1\simeq S^1$, we may
view $PSL(2,\R)$ as a subgroup of $\mbox{Homeo}_+(S^1)$ and
$\widetilde{PSL}(2,\R)$ as a subgroup of
$\widetilde{\mbox{Homeo}}_+(S^1)$. Here $\widetilde{PSL}(2,\R)$ is
viewed as the group of equivalence classes of paths in $PSL(2,\R)$
starting at the identity. Also let $\widetilde{PSL}(2,\Z)\subset
\widetilde{PSL}(2,\R)$ be the equivalence classes of paths in
$PSL(2,\R)$ starting at $id$ and ending at an element in
$PSL(2,\Z)$. We now restrict $\Psi$ to $\widetilde{PSL}(2,\Z)$,
which is isomorphic to $Aut(S,\bdry S)$, where $S$ is the
once-punctured torus. Recall that any $\sigma\in B_3\cong
Aut(S,\bdry S)$ can uniquely be written as
$(\sigma_1\sigma_2\sigma_1)^k w$, where $w$ is a product of
$\sigma_1, \sigma_2^{-1}$ or $\sigma_1^{-1},\sigma_2$. Since $\Psi$
and $-4~rot$ agree on all powers of $(\sigma_1\sigma_2\sigma_1)^4$,
they differ by a bounded amount on $\widetilde{PSL}(2,\Z)$. Hence
their coboundaries $\delta\Psi$ and $\delta(-4~rot)$ represent the
same element in the bounded cohomology group $H^2_b(PSL(2,\Z);\Z)$.

Next, given the Rademacher function $\Phi: PSL(2,\Z)\rightarrow \Z$,
consider its coboundary $\delta\Phi$. Although $\delta\Phi$ is zero
in the ordinary group cohomology $H^2(PSL(2,\Z);\Z) =\Z/6\Z$, it is
nevertheless a nontrivial element in the bounded cohomology group
$H^2_b(PSL(2,\Z);\Z)$:  First observe that $\Phi$ is not a bounded
$1$-cochain.  Moreover, since $PSL(2,\Z)\cong \Z/2\Z * \Z/3\Z$,
there is no nonzero homomorphism $PSL(2,\Z)\rightarrow \Z$. (Observe
that there is no nonzero homomorphism from $\Z/m\Z$ to $\Z$, when
$m$ is a positive integer.) Hence $\Phi$ is not a bounded distance
from any homeomorphism and therefore represents a nontrivial element
in $H^2_b(PSL(2,\Z);\Z)$.

Consider the following diagram --- keep in mind that we need to
distinguish among similar-looking groups $PSL(2,\Z)$, $PSL(2,\R)$,
etc.:

$$\begin{CD}
0 @> >> \Z @> >> Aut(S,\bdry S)=\widetilde{PSL}(2,\Z) @>  >> PSL(2,\Z)   @>   >>0  \\
&& @V VV @V  VV @V  VV \\
0 @> >> \Z @> >> \widetilde{PSL}(2,\R) @>  >>  PSL(2,\R)   @>   >> 0 \\
\end{CD}   $$

Theorem~\ref{lk} implies the following:

\begin{cor}
$\delta \Phi=-12\delta(\rot)$ as 2-cochains on $PSL(2,\Z)$.
\end{cor}

In other words, two seemingly different quasi-morphisms
--- the translation number for $\widetilde{\mbox{Homeo}}_+(S^1)$
and the Rademacher function --- have essentially the same
coboundary. Hence, we can keep track of the value of one
quasi-morphism through the value of the other quasi-morphism,
although the functions are far from identical.

\subsection{The Maslov index}\label{maslov}
In this subsection, we define the {\em Maslov index}. There are
various definitions of the Maslov index in the literature, and our
$\mu(\gamma,\Lambda_0)$ is identical to that of  Robbin and Salamon
in \cite{RS}.

Consider the symplectic vector space $(\R^{2n}=\R^n\times
\R^n,\omega)$, with coordinates $\mathbf{x}=(x_1,\dots,x_n)$ for the
first $\R^n$ and $\mathbf{y}=(y_1,\dots,y_n)$ for the second $\R^n$,
and symplectic form $\omega=\sum_{i=1}^n dx_i \wedge dy_i$. Consider
the Lagrangian subspace $\Lambda_0=\{\mathbf{y}=0\}$. Let
$\mathcal{L}$ be the {\em Lagrangian Grassmannian} of
$(\R^{2n},\omega)$, i.e., the set of Lagrangian subspaces of
$\R^{2n}$. Also let $\mathcal{L}_{\Lambda_0}$ be the {\em Maslov
cycle} of $\Lambda_0$, namely the set of $\Lambda\in \mathcal{L}$
with $\Lambda_0\cap \Lambda\not=\{0\}$.

If $\Lambda\in \mathcal{L}$ is transverse to $\Lambda_0$, then there
exists an element of the symplectic group $Sp(2n,\R)$ which sends
$\Lambda_0$ to itself and $\Lambda$ to
$\Lambda_0'=\{\mathbf{x}=0\}$. [Proof: If $v_1,\dots,v_n$ is a basis
for $\Lambda_0$ and $w_1,\dots,w_n$ is a basis for $\Lambda$, then,
with respect to $v_1,\dots,v_n,w_1,\dots,w_n$, the symplectic form
can be written as $\begin{pmatrix} 0 & A\\ -A^T & 0
\end{pmatrix}$, for some nonsingular $n\times n$ matrix $A$.  Since
we are allowed to change bases of $\Lambda_0$ and $\Lambda$
(separately), we are looking to solve:
$$\begin{pmatrix} B^T & 0 \\ 0 & C^T \end{pmatrix}
\begin{pmatrix} 0 & A\\ -A^T & 0 \end{pmatrix}
\begin{pmatrix} B & 0 \\ 0 & C\end{pmatrix}
= \begin{pmatrix} 0 & B^T AC \\ -C^T A^T B & 0\end{pmatrix}
=\begin{pmatrix} 0 & I \\ -I & 0\end{pmatrix}.$$ Such $B,C$ can
easily be found.]

Now consider the neighborhood $U=\{\mathbf{y}=A\mathbf{x}~|~ A
\mbox{ symmetric $n\times n$ matrix} \}$ of $\Lambda_0\in
\mathcal{L}$. [Proof: Any $n$-plane which is sufficiently close to
$\mathbf{y}=0$ is graphical of form $\mathbf{y}=A\mathbf{x}$.  We
can check that the Lagrangian condition implies that $A=A^T$.] We
observe that $U$ depends on the choice of complementary Lagrangian
subspaces $\Lambda_0$ and $\Lambda_0'$, and will often be written as
$U(\Lambda_0,\Lambda_0')$.  It is easy to see that $U$ is
contractible.

A useful operation which allows us to cover all of $\mathcal{L}$
with open sets of type $U(\Lambda_0,\Lambda_0')$, is the {\em
symplectic shear} $\begin{pmatrix} I & A\\ 0 & I \end{pmatrix}\in
Sp(2n,\R)$, where $A$ is a symmetric $n\times n$ matrix.  The shear
sends $\Lambda_0$ to itself and $\{\mathbf{x}=0\}$ to
$\{\mathbf{x}=A\mathbf{y}\}$.  If $A$ is invertible, then the
Lagrangian subspace $\{(A\mathbf{y}, \mathbf{y})~|~\mathbf{y}\in
\R^n\}=\{(\mathbf{x}, A^{-1}\mathbf{x})~|~\mathbf{x}\in \R^n\}$ is
in $U(\Lambda_0,\{\mathbf{x}=0\})$. By ranging $\Lambda_0'$ over all
the Lagrangian subspaces transverse to $\Lambda_0$, the collection
of such $U(\Lambda_0,\Lambda_0')$ covers $\mathcal{L}$.

We can now define the {\em Maslov index} of a path
$\gamma:[0,1]\rightarrow \mathcal{L}$, with respect to a fixed
Lagrangian $\Lambda_0$.  Subdivide $[0,1]$ into
$0=t_0<t_1<\dots<t_k=1$, so that each $\gamma|_{[t_i,t_{i+1}]}$ lies
in some $U(\Lambda_0,\Lambda_0')$.  Suppose $\gamma(t_i)$ is given
by $\{\mathbf{y}=A(t_i) \mathbf{x}\}$ and $\gamma(t_{i+1})$ by
$\{\mathbf{y}=A(t_{i+1}) \mathbf{x}\}$.  Then let
\begin{equation}\label{define-Maslov}
\mu(\gamma|_{[t_{i},t_{i+1}]}, \Lambda_0)={1\over 2} sign(A(t_{i+1}))
-{1\over 2}sign(A(t_{i})).
\end{equation}
Here $sign$ denotes the signature of the symmetric matrix. (The
signature of a symmetric bilinear form is the dimension of the
maximal positive definite subspace minus the dimension of the
maximal negative definite subspace.) We then define
\begin{equation}
\mu(\gamma,\Lambda_0)\stackrel{def}=\sum_{i=0}^{k-1}
\mu(\gamma|_{[t_i,t_{i+1}]},\Lambda_0).
\end{equation}
By Theorem~2.3 of \cite{RS}, this $\mu$ is well-defined, invariant
under homotopies fixing endpoints, and is natural, i.e.,
$\mu(\Psi(\gamma), \Psi(\Lambda))=\mu(\gamma,\Lambda)$, where
$\Psi\in Sp(2n,\R)$.  Moreover, if $\gamma$ does not intersect the
Maslov cycle $\mathcal{L}_{\Lambda_0}$, then
$\mu(\gamma,\Lambda_0)=0$. In the special case that $\gamma$ is a
closed loop, $\mu(\gamma, \Lambda_0)$ is independent of the choice
of $\Lambda_0$.

Next, given $L_1,L_2,L_3\in \mathcal{L}$, we define the {\em ternary
index} $I(L_1,L_2,L_3)$.  Consider the symmetric bilinear form $Q$
on $(L_1+L_2)\cap L_3$ defined by $Q(v,w)=\omega(v_2,w)$, where
$v\in (L_1+L_2)\cap L_3$ is written as $v_1+v_2$, with $v_1\in L_1$,
$v_2\in L_2$. Then $I(L_1,L_2,L_3)$ is the signature of $Q$.

It is not difficult to see that $I(L_1,L_2,L_3)$ has the following
equivalent definition:  Consider the subspace $V\subset L_1\oplus
L_2\oplus L_3$, consisting of triples $(v_1,v_2,v_3)$, $v_i\in L_i$,
with $v_1+v_2+v_3=0$.  Define the quadratic form $Q':V\times
V\rightarrow \R$ by:
\begin{eqnarray*}
Q'((v_1,v_2,v_3),(w_1,w_2,w_3)) &=& \omega(v_1,w_3)=\omega(v_2,w_1)=\omega(v_3,w_2)\\
&=& -\omega(v_3,w_1) =-\omega(v_1,w_2) = -\omega(v_2,w_3).
\end{eqnarray*}
Then $I(L_1,L_2,L_3)$ is also the signature of $Q'$.

Now, given $L_1,L_2,L_3\in \mathcal{L}$, let $\gamma_{12}$ be a path
in $\mathcal{L}$ from $L_1$ to $L_2$, $\gamma_{23}$ be a path in
$\mathcal{L}$ from $L_2$ to $L_3$, and let $\gamma_{13}$ be the
concatenation $\gamma_{12}\gamma_{23}$.  Also let
$\gamma_{ij}=\gamma_{ji}^{-1}$.  We then have the following:

\begin{lemma}
$I(L_1,L_2,L_3)= 2(\mu(\gamma_{12}, L_1)+\mu(\gamma_{23},L_2)+\mu(\gamma_{31},L_3)).$
\end{lemma}

\begin{proof}
Suppose $L_1$, $L_2$, $L_3$ are mutually transverse.  Without loss
of generality, we may take $L_1=\{\mathbf{y}=0\}$,
$L_2=\{\mathbf{x}=0\}$, and $L_3=\{\mathbf{y}=A\mathbf{x}\}$, where
$A$ is symmetric and nonsingular.

Since the right-hand side of the equation in the lemma does not
depend on the choice of paths, provided the endpoints remain the
same, there is no loss of generality in proving the lemma for a
convenient choice of paths. (As remarked earlier, if $\gamma$ is a
loop, then $\mu(\gamma,L)$ does not depend on the choice of $L\in
\mathcal{L}$.)  Define $\gamma_{13}:[0,1]\rightarrow \mathcal{L}$ by
$t\mapsto \{\mathbf{y}=tA\mathbf{x}\}$,
$\gamma_{23}:[0,1]\rightarrow \mathcal{L}$ by $t\mapsto
\{\mathbf{x}=tA^{-1}\mathbf{y}\}$, and
$\gamma_{12}=\gamma_{13}\gamma_{32}$. One easily computes from
Equation~\ref{define-Maslov} that
\begin{eqnarray*}
\mu(\gamma_{12},L_1)&=&{1\over 2}sign(A),\\
\mu(\gamma_{23},L_2)&=&-{1\over 2}sign(A^{-1})=-{1\over 2}sign(A),\\
\mu(\gamma_{31},L_3)&=&-{1\over 2}sign(A).
\end{eqnarray*}
Therefore, the right-hand side of the equation in the lemma equals $-sign(A)$.

On the other hand, $(L_1+L_2)\cap L_3=L_3$ in our case, and
$$Q((\mathbf{x},A\mathbf{x}),(\mathbf{x}',A\mathbf{x}'))=\omega((0,A\mathbf{x}),
(\mathbf{x}',0))=-\mathbf{x}^T A^T \mathbf{x}'.$$ Thus,
$I(L_1,L_2,L_3)=-sign(A)$ as well.

The general case is more involved, and we only sketch the idea.
First, we normalize $L_1,L_2,L_3$ (this requires some work), and
then use the {\em additivity formula} from \cite{RS}. Let $V'$ be a
Lagrangian of standard symplectic $\R^{2n'}$ and $V''$ be a
Lagrangian of $\R^{2n''}$. Also let $\gamma'$ be a path in
$\mathcal{L}(\R^{2n'})$ and $\gamma''$ be a path in
$\mathcal{L}(\R^{2n''})$.  Then the additivity formula says the
following:
\begin{equation}
\mu(\gamma'\oplus \gamma'',V'\oplus V'')=
\mu(\gamma',V')+\mu(\gamma'',V'').
\end{equation}
We can then reduce to the above calculation where $L_1,L_2,L_3$ are
mutually transverse.
\end{proof}

\subsection{The signature}\label{thesignature}
To generalize the results we obtained for the punctured torus to the
braid group $B_n$ on $n$ strands, we use the {\em signature}. The
signature has the useful property of remaining invariant under
stabilization, whereas the linking number increases under
stabilization.

Define the {\em signature function} on $B_n$ as follows: Given a
braid $\alpha$, let $\widehat\alpha$ denote the braid closure inside
$S^3$, and $sign(\widehat\alpha)$ be the signature of the
(symmetrized) Seifert pairing. It is not difficult to see that the
signature is a quasi-morphism on $B_n$.

We can view $B_{2n+1}$ as the hyperelliptic mapping class group of a
once-punctured surface $\Sigma_n$, where $n$ is the genus, and
$B_{2n}$ as the hyperelliptic mapping class group of a
twice-punctured surface $\Sigma_{n-1}$, where $n-1$ is the genus.
(Here by the genus we mean the genus of the closed surface obtained
by adding disks.)  For $B_{2n+1}$, let $\mathcal{B}_{-1}$ be the map
$B_{2n+1}\rightarrow Sp(2n,\R)$, which is the action on the
symplectic vector space  $H_1(\Sigma_n,\bdry\Sigma_n;\R)$ (with
symplectic form the intersection pairing).  For $B_{2n}$, the
intersection pairing is degenerate, so we take the standard
embedding of $B_{2n}$ into $B_{2n+1}$ by adding a trivial strand;
from now on we assume that our braid groups have an odd number of
stands.  The strange notation $\mathcal{B}_{-1}$ comes from the fact
that the homology representation is the Burau representation
specialized at $-1$. Note that for $n=1$ it is the same as the map
$B_{3}\rightarrow SL(2,\Z)$ which appeared in Section 2.

Next we define the {\em Meyer cocycle} $Meyer(g_1,g_2)$, where
$g_1,g_2\in Sp(2n,\R)$.  (Here we are thinking of
$H_1(\Sigma_n,\bdry\Sigma_n;\R)$ as a symplectic vector space of
dimension $2n$.) Consider the symplectic vector space
$(\R^{2n}\times \R^{2n}, \omega\oplus -\omega)$.  Let $\widetilde
g_1$ be a path in $Sp(2n,\R)$ from $id$ to $g_1$ and $\widetilde
g_2$ be a path from $id$ to $g_2$.  Also let
$\widetilde{g_1g_2}(t)=\widetilde g_1(t)\widetilde g_2(t)$,
$t\in[0,1]$; this is homotopic to the path $\widetilde g_1(t)$,
followed by the path $g_1\cdot(\widetilde g_2(t))$.  Now let
$Graph(h)$ be the graph of $h\in Sp(2n,\R)$, i.e., it is the
Lagrangian of $\R^{2n}\oplus \R^{2n}$ consisting of vectors
$(v,h(v))$; if $\widetilde h$ is a path in $Sp(2n,\R)$,
$Graph(\widetilde h)$ is a path of Lagrangians $Graph(\widetilde
h(t))$, $t\in[0,1]$.  Then we set:
\begin{eqnarray*}
Meyer(g_1,g_2)&\stackrel{def}=& I(Graph(id),Graph(g_1),Graph(g_1g_2))\\
&=& 2(\mu(Graph(\widetilde g_1),Graph(id))+\mu(Graph(\widetilde g_2),Graph(id))\\
& & ~~~~ -\mu(Graph(\widetilde{g_1g_2}),Graph(id))).
\end{eqnarray*}
[It is not hard to verify that if $\gamma:[0,1]\rightarrow
\mathcal{L}$ is a path of Lagrangians, then
$\mu(\gamma,\gamma(0))=-\mu(\gamma^{-1},\gamma^{-1}(0))$, where
$\gamma^{-1}(t)=\gamma(1-t)$. We also used the fact that
$\mu(\gamma,\Lambda)=\mu(\Psi(\gamma),\Psi(\Lambda))$]

In \cite{GG1} it is proven that:
\begin{equation}\label{signature}
sign(\widehat{\alpha\beta})= sign(\widehat\alpha)
+sign(\widehat\beta) -
Meyer(\mathcal{B}_{-1}(\alpha),\mathcal{B}_{-1}(\beta)).
\end{equation}
There are two quasi-morphisms whose coboundary is the Meyer cocycle:
on $B_{2n+1}$ there is the {\em signature}, and on
$\widetilde{Sp}(2n,\R)$ there is the {\em Maslov index}.  Here
$\widetilde{Sp}(2n,\R)$ is the universal cover of $Sp(2n,\R)$.  In
order to relate the two, we first observe the following:

\begin{lemma}
The homomorphism $\mathcal{B}_{-1}:B_{2n+1}\rightarrow Sp(2n,\R)$
can be lifted to a homomorphism
$\widetilde{\mathcal{B}}_{-1}:B_{2n+1}\rightarrow
\widetilde{Sp}(2n,\R)$.
\end{lemma}

\begin{proof}
Let $C_1,\dots,C_{2n}$ be oriented nonseparating closed curves on
$\Sigma_{n}$ so that (i) $\sigma_i$ corresponds to a positive Dehn
twist about $C_i$ and (ii) the intersection pairing
$\omega(C_i,C_j)$ is $\delta_{i+1,j}- \delta_{i-1,j}$. (Here
$\delta_{i,j}$ is $0$ if $i\not=j$ and $1$ if $i=j$.)

Let $A_1$ be the $2\times 2$ matrix $\begin{pmatrix} 1 & 1\\ 0 &
1\end{pmatrix}$ and let $A_2$ be the $3\times 3$ matrix
$\begin{pmatrix} 1 & 0 & 0\\ -1 & 1 & 1\\ 0 & 0 & 1 \end{pmatrix}$.
Then $\mathcal{B}_{-1}(\sigma_1)=diag(A_1,1,\dots,1)$. (By this we
mean the matrix which has the given entries along the diagonal and
has zeros everywhere else.) We also have
$\mathcal{B}_{-1}(\sigma_2)=diag(A_2,1,\dots,1), \dots,
\mathcal{B}_{-1}(\sigma_{2n-1})=diag(1,\dots,1,A_2)$, and
$\mathcal{B}_{-1}(\sigma_{2n})=diag(1,\dots,1,A_1)$.

To lift to $\widetilde{Sp}(2n,\R)$, we replace $A_1$ by
$A_1(t)=\begin{pmatrix} 1 & t\\ 0 & 1\end{pmatrix}$ and $A_2$ by
$A_2(t)=\begin{pmatrix} 1 & 0 & 0\\ -t & 1 & t\\ 0 & 0 & 1
\end{pmatrix}$, where $t\in[0,1]$.  (We do this for all the
$\mathcal{B}_{-1}(\sigma_i)$.)  To verify that this indeed gives a
lift $\widetilde{\mathcal{B}}_{-1}: B_n\rightarrow
\widetilde{Sp}(2n,\R)$, we need to check the braid relations.

If $|i-j|\geq 2$, then $\sigma_i\sigma_j=\sigma_j\sigma_i$, and to
verify that
$\widetilde{\mathcal{B}}_{-1}(\sigma_i)\widetilde{\mathcal{B}}_{-1}(\sigma_j)
=\widetilde{\mathcal{B}}_{-1}(\sigma_j)\widetilde{\mathcal{B}}_{-1}(\sigma_i)$
it suffices to check that $diag(A_2(t),1,1)$ and $diag(1,1,A_2(t))$
commute. This is an easy calculation. (The cases $i=1, j=3$ and $i=
2n-2, j=n$, when $diag(A_1(t),1,1)$ and $diag(1,1,A_1(t)) $ are
involved, are easier.)

We also verify that
$$\widetilde{\mathcal{B}}_{-1}(\sigma_i)\widetilde{\mathcal{B}}_{-1}(\sigma_{i+1})
\widetilde{\mathcal{B}}_{-1}(\sigma_i)
=\widetilde{\mathcal{B}}_{-1}(\sigma_{i+1})
\widetilde{\mathcal{B}}_{-1}(\sigma_{i})
\widetilde{\mathcal{B}}_{-1}(\sigma_{i+1}).$$ Computing both sides,
we require:
\begin{equation}\label{matrices}
\begin{pmatrix} 1 & 0 & 0 & 0\\ -2t+t^3 & 1-t^2 & 2t-t^3 & t^2\\
t^2 & -t & 1-t^2 & t\\
0 & 0 & 0 & 1\end{pmatrix}\mbox{ and }
\begin{pmatrix} 1 & 0 & 0 & 0 \\ -t & 1-t^2 & t & t^2\\
t^2 & -2t+t^3 & 1-t^2 & 2t-t^3\\
0 & 0 & 0 & 1\end{pmatrix}
\end{equation}
to be homotopic as paths. Since $-t$ and $-2t+t^3$ are both negative
for $t\in(0,1]$, we can take $a(s,t)=(1-s)(-2t+t^3)+s(-t)$ and
$b(s,t)={(-t)(-2t+t^3)\over a(s,t)}$.  Now
$$\begin{pmatrix} 1 & 0 & 0 & 0 \\ a(s,t) & 1-t^2 & -a(s,t) & t^2\\
t^2 & b(s,t) & 1-t^2 & -b(s,t)\\
0 & 0 & 0 & 1\end{pmatrix}$$ is a homotopy of paths in $Sp(2n,\R)$
which takes the left-hand matrix in Equation~\ref{matrices} to the
right-hand one in Equation~\ref{matrices}.
\end{proof}

The following is a relatively simple computation, once the
definitions are sorted out:

\begin{lemma}\label{onedehn}
$\mu(Graph(\widetilde h),Graph(id))={1\over 2}$ if $h\in Sp(2n,\R)$
is a positive Dehn twist about a nonseparating curve.
\end{lemma}

\begin{proof}
We first reduce to the case where $n=1$, the symplectic form on
$\R^2\oplus \R^2$ (with coordinates $(x,y)=((x_1,x_2),(y_1,y_2))$)
is $\omega=dx_1\wedge dx_2-dy_1\wedge dy_2$, and
$h=\begin{pmatrix}
1 & 1\\0 & 1
\end{pmatrix}$.  Indeed, $Graph(id)$ is the set
$\{(x,x)~|~x\in\R^{2n}\}=\{(gx,gx)~|~ x\in \R^{2n}\}$ and $Graph(h)$
is the set $\{(x,hx)~|~x\in \R^{2n}\}=\{(gx,hgx)~|~ x\in \R^{2n}\}$,
if $g$ is a nonsingular $2n\times 2n$ matrix.  Now, apply
$(g^{-1},g^{-1})$ to both $Graph(id)$ and $Graph(h)$, where $g\in
Sp(2n,\R)$.  This gives us $Graph(id)$ and $Graph(g^{-1}hg)$. Hence,
by conjugating, we may assume that $h$ is as above, since $\mu$ is
invariant under the action of the symplectic group.

The graph of $id$, which we write as $L_0$, is
$\R\{v_1=(1,0,1,0),v_2=(0,1,0,1)\}$.  A complementary Lagrangian
subspace to $L_0$ is $L_0'=\R\{w_1=(0,0,0,-1), w_2=(-1,0,0,0)\}$.
(Here $\omega(v_i,w_j)=\delta_{i,j}$.) The graph of $h$ is spanned
by $v_1=(1,0,1,0)$ and $v_1+v_2+w_2=(0,1,1,1)$, or, equivalently, by
$v_1$ and $v_2+w_2$.  Hence $$\mu(Graph(\widetilde
h),Graph(id))={1\over 2} sign\begin{pmatrix} 0 & 0 \\ 0 & 1
\end{pmatrix}={1\over 2}.$$ This proves the lemma.
\end{proof}

We now state the main theorem of this section:

\begin{thm}\label{sign+maslov}
Let $\gamma$ be an element of $B_{2n+1}$, or equivalently, an
element of a hyperelliptic mapping class group $Hyp
Aut(\Sigma_n,\bdry \Sigma_n)$. Then $$sign(\widehat{\gamma})=
-lk(\gamma)+2\mu(Graph(\widetilde{\mathcal{B}}_{-1}(\gamma)),Graph(id)).$$
\end{thm}

Since $lk(\gamma)\geq 0$ if $\gamma$ is quasipositive, we have the
following:

\begin{cor}\label{cor: restriction}
If $sign(\widehat\gamma) >
2\mu(Graph(\widetilde{\mathcal{B}}_{-1}(\gamma)),Graph(id)),$ then
$\gamma$ cannot be quasipositive in $B_{2n+1}$. Equivalently,
$\gamma$ is not a product of positive Dehn twists in $Hyp Aut(
\Sigma_n, \bdry\Sigma_n)$.

\end{cor}

\begin{proof}
Let $\gamma$ be an element of $B_{2n+1}$.  Then $\gamma$ can be
written as $\gamma=\gamma_1\dots\gamma_k$, where $\gamma_i$ are all
conjugates of a standard half-twist or its inverse. Let
$g_i=\mathcal{B}_{-1}(\gamma_i)\in Sp(2n,\R)$, and let $\widetilde
g_i=\widetilde{\mathcal{B}}_{-1}(\gamma_i)\in \widetilde{Sp}(2n,\R)$
be a path from $id$ to $g_i$ in $Sp(2n,\R)$.  By repeatedly using
Equation~\ref{signature} and observing that
$sign(\widehat\gamma_i)=0$ (since the Seifert surface is a disk), we
have:
\begin{eqnarray*}
sign(\widehat{\gamma_1\dots\gamma_k}) &=& sign(\widehat\gamma_1)
+sign(\widehat{\gamma_2\dots\gamma_k}) -
Meyer (\mathcal{B}_{-1}(\gamma_1),\mathcal{B}_{-1}(\gamma_2\dots\gamma_k))\\
&=&sign(\widehat{\gamma_2\dots\gamma_k}) -
Meyer (\mathcal{B}_{-1}(\gamma_1),\mathcal{B}_{-1}(\gamma_2\dots\gamma_k))\\
&=& - \sum_{i=1}^{k-1} Meyer(g_i, g_{i+1}\dots g_k).
\end{eqnarray*}
Next, applying Lemma~\ref{onedehn}, we have:
\begin{eqnarray*}
sign(\widehat{\gamma_1\dots\gamma_k})&=& -2\sum_{i=1}^{k-1}
\left\{ \mu(Graph(\widetilde g_i),Graph(id))+
\mu(Graph(\widetilde{g_{i+1}\dots g_k}),Graph(id))\right.\\
& & \hspace{.5in} \left. - \mu(Graph(\widetilde{g_{i}\dots g_k}),Graph(id))\right\}\\
&=& -2\sum_{i=1}^{k-1} \left\{ \pm {1\over 2}+
\mu(Graph(\widetilde{g_{i+1}\dots g_k}),Graph(id))\right.\\
& & \hspace{.5in} \left. - \mu(Graph(\widetilde{g_{i}\dots g_k}),Graph(id))\right\}\\
&=& -lk(\gamma)+2\mu(Graph(\widetilde{g_1\dots g_k}),Graph(id)).
\end{eqnarray*}
Here we have $\pm {1\over 2}$ depending on whether we have a
positive or negative Dehn twist.

If $\gamma$ is quasipositive, then $lk(\gamma)\geq 0$. Hence
$sign(\gamma)\leq 2\mu (Graph(\widetilde{g_1\dots g_k}),Graph(id))$.
\end{proof}

\s\n {\bf Remark.}  In \cite{GG1}, Gambaudo and Ghys prove that, for
the ``generic element'' $\gamma\in B_3$,
$$sign(\widehat\gamma) +{2\over 3} lk(\gamma) =-{1\over 3} \Phi(\mathcal{B}_{-1}(\gamma)).$$
For example, if $\gamma$ is generic if it is pseudo-Anosov.
Combining with Theorem~\ref{lk}, we have:
$$sign(\widehat\gamma)=-lk(\gamma) + 4~\rot(\gamma),$$
for such $\gamma$. This is consistent with
Theorem~\ref{sign+maslov}.

\s\n {\bf Remark.}  Gambaudo and Ghys also have a formula analogous
to Equation~\ref{signature} for the $\omega$-signatures.  Presumably
our Theorem~\ref{sign+maslov} can be generalized to
$\omega$-signatures as well.

\section{Characterization of (weak) symplectic fillability}
\label{section: weak fillability}

In this section we prove Theorem~\ref{cgeq1}. The starting point is
the following special case of a theorem of Roberts~\cite{Ro1, Ro2},
generalizing work of Hatcher~\cite{Ha}.

\begin{thm}[Roberts] \label{ourversion} Assume the surface $S$ has one
boundary component and $h$ is a diffeomorphism that restricts to the
identity on the boundary. If $h$ is isotopic to a pseudo-Anosov
homeomorphism $\psi$ and the fractional Dehn twist coefficient of
$h$ is $c$, then $M=(S,h)$ carries a taut foliation transverse to
the binding if $c\ge1$.
\end{thm}

This theorem (not stated in this form by Roberts) follows from a
more general result of Roberts which is stated below as
Theorem~\ref{Roberts}.  To explain how Theorem~\ref{ourversion}
follows from Theorem~\ref{Roberts}, we start by comparing the
notation and coordinates used by Roberts to our own.

Let $S$ be a hyperbolic surface with one boundary component, $\psi$
be the pseudo-Anosov representative of $h$, and $c$ be the
fractional Dehn twist coefficient. Denote by $N$ the mapping torus
of $\psi$, i.e., $N \stackrel{def} = (S\times [0,1])/(x,
1)\sim(\psi(x),0)$ for $x\in S$.

Roberts gives an oriented identification $\bdry N\simeq \R^2/\Z^2$
by choosing closed curves $\lambda, \mu$ so that $\lambda$ has slope
$0$ and $\mu$ has slope $\infty$. (See Section 3 of \cite{Ro2}.)
Here we choose orientations to agree with the usual conventions for
a knot complement. We will now describe the curves $\lambda$ and
$\mu$.  Let $\lambda=\bdry (S\times\{0\})$.  Define $\gamma$ to be
one component of the suspension of the periodic points of
$\psi|_{\bdry S}$.  If there are $n$ prongs, then there are $2n$
periodic points, $n$ of which are attracting and $n$ of which are
expanding. Observe that the geometric intersection number
$\#(\gamma\cap (S\times\{0\}))$ divides $n$ and equals $n$ if the
suspension (of only the attracting points) is connected. Now we
define $\mu$ to be the essential closed curve on $\bdry N$ which has
the minimal $\#(\mu\cap \gamma)$ amongst all closed curves on $\bdry
N$ which form an integral basis of $H_1(\bdry N;\Z)$ with $\lambda$.
The choice of $\mu$ is not unique if $\#(\gamma\cap\lambda)=2$;
there are two choices which minimize $\#(\mu\cap\gamma)$.  In that
case we choose $\mu$ so that $\mbox{slope}(\gamma)=+2$.

We now state Theorem 4.7 of \cite{Ro2}:

\begin{thm}[Roberts] \label{Roberts}
Suppose $S$ has one boundary component, $\psi$ is a pseudo-Anosov
map and $N=(S\times[0,1])/(x,1)\sim(\psi(x),0)$.  Then one of the
following holds:
\be
\item $\gamma$ has slope infinity and $N$ contains taut foliations
realizing all boundary slopes in $(-\infty, \infty)$.

\item $\gamma$ has positive slope and $N$ contains taut foliations
realizing all boundary slopes in $(-\infty,1)$.

\item $\gamma$ has negative slope and $N$ contains taut foliations
realizing all boundary slopes in $(-1,\infty)$. \ee Here the slope
is measured with respect to the identification $\bdry N\simeq
\R^2/\Z^2$ given by the basis $(\lambda, \mu)$ defined above, and
``realizing'' a boundary slope means the restriction of the taut
foliation to $\bdry N$ is a linear foliation with the given boundary
slope.
\end{thm}

\begin{proof}
[Proof that Theorem~\ref{Roberts} implies Theorem~\ref{ourversion}]

Suppose $h \in  Aut(S,\bdry S)$, $\psi$ is its pseudo-Anosov
representative, and $c={p\over q}$ is the corresponding fractional
Dehn twist coefficient.  Assume $p,q$ are relatively prime positive
integers. If the closed manifold $M=(S,h)$ is obtained by Dehn
filling $N$ along the closed curve $\nu$ on $\bdry N$, then $\gamma=
p \lambda + q \nu$ in $H_1(\bdry N;\Z)$. We also have $\mu = \nu +
k\lambda$, where $k$ is an integer chosen to minimize $|\gamma\cdot
\mu| = |(p\lambda +q\nu)\cdot (\nu+k\lambda)|=|p-kq|$. When there is
a tie, i.e., both $p-kq=\pm {p\over 2}$ are possible, the tie is
broken by choosing ${p\over 2}$. In the cases below, the slope will
be computed relative to the basis $(\lambda, \mu)$. We compute that
slope$(\nu)={(\lambda \cdot \nu) \over (\nu\cdot \mu)} =-{1\over k}$
and slope$(\gamma)={q\over p-kq}$.

\s\n (1) $c={p\over q}$ is an integer $\ge 1 $.  It follows that
$p-kq=0$ and $k\ge 1$. Therefore, slope$(\gamma)=\infty$ and
slope$(\nu)=-{1 \over k} \in[-1,0) \subset (-\infty, \infty)$.

\s\n (2) ${p\over q} >1$ is not an integer and $k$ satisfies $0<
p-kq \le {p\over 2}$.  It follows that slope$(\gamma)>0$, $k\ge 1$,
and therefore slope$(\nu)= -{1\over k} \in[-1,0) \subset (-\infty,
1)$.

\s\n (3) ${p\over q} > 1$ is not an integer and $k$ satisfies
$-{p\over 2} < p-kq < 0 $.  It follows that slope$(\gamma)<0$, $k
\ge 2$, and therefore slope$(\nu)= -{1\over k} \in[-{1\over 2}, 0)
\subset (-1, \infty)$.

\s\n Thus, for all $c\ge1$, a taut foliation of $N$ can be
constructed with boundary slope equal to the slope of the meridian
of the solid torus that extends $N$ to $M$. By extending the leaves
by meridian disks, we can construct a taut foliation of $M$
transverse to the binding.
\end{proof}

\s If $M=(S,h)$ and $c\ge 1$, let $\mathcal F$ be a taut foliation
furnished by Theorem~\ref{ourversion}.  By the work of
Eliashberg-Thurston \cite{ET}, any taut foliation admits a
$C^0$-small perturbation into a universally tight and (weakly)
symplectically fillable contact structure.  We denote a perturbation
of $\mathcal{F}$ by $\xi_{\mathcal{F}}$.  (Note that, a priori, two
perturbations of $\mathcal{F}$ may not even be isotopic.)  We will
denote by $(S,h)$ the contact structure corresponding to the open
book, which is also denoted $(S,h)$.

\begin{thm}\label{OB=ET}
If $c\ge 1$, then the contact structure $(S,h)$ is isotopic to
$\xi_{\mathcal{F}}$ for some taut foliation $\mathcal{F}$.
\end{thm}

The proof of Theorem~\ref{OB=ET} will occupy the rest of the
section. We first claim the following:

\begin{lemma}\label{vf}
Let $\psi$ be a pseudo-Anosov representative of $h \in Aut(S,\bdry
S)$ with $c \ge 1$ and let $N  = (S\times [0,1])/(x,
1)\sim(\psi(x),0)$. Then there exists a nonsingular vector field $X$
on $N$ with the following properties:
\be
\item $X$ is tangent to $\bdry N$.
\item $X$ is positively transverse to $S\times\{t\}$ for all $t\in[0,1]$.
\item There exists a transversely oriented taut foliation $\mathcal{F}$
on $N$ which is positively transverse to $X$. Moreover,
$\mathcal{F}$ can be chosen so that $\mathcal{F}\cap T(\bdry N)$ is
a nonsingular foliation on $\bdry N$ which is foliated by circles of
slope $-{1\over k}$, where $k$ is the positive integer as described
in the proof of Theorem~\ref{Roberts}. \ee
\end{lemma}

Here we are using slope convention used in Theorem~\ref{Roberts}.

\begin{proof}
This follows from analyzing Roberts' construction (cf. Section 2 of
\cite{Ro1}) of the foliation $\mathcal{F}$ in
Theorem~\ref{ourversion} and noting that it can be performed in a
manner compatible with $X$. Roberts constructs a set
$\alpha_1,\dots,\alpha_m$ of properly embedded oriented arcs in $S$
with the following properties. Consider $D_i=\alpha_i\times
[{i-1\over m},{i\over m}]$, where $D_i$ is oriented so that
$({\bdry\over \bdry t},\dot{\alpha})$ form an oriented basis for
$TD_i$. (Here $t$ is the coordinate for $[0,1]$.) Also write
$S_t\stackrel{def} =S\times\{t\}$. This gives a spine
$$\Sigma=(\cup_{i=1}^m S_{i/m})\cup (\cup_{i=1}^m D_i),$$
which can be modified into a branched surface $\mathcal{B}$ by
isotoping $D_i$ and smoothing the neighborhood of each intersection
$\alpha_i\times \{ {i\over m}\}$ between $S_{i/m}$ and $D_i$ into a
branch locus, so that, near the branch locus, a vector field which
is positively transverse to $S_{i/m}$ (we may take ${\bdry\over
\bdry t}$ here) is also positively transverse to the new $D_i$.  The
same can be done for each intersection $\alpha_i\times \{ {i-1\over
m} \}$ between $S_{(i-1)/m}$ and $D_{i}$.

On each $S\times[{i-1\over m},{i\over m}]$, start with ${\bdry\over
\bdry t}$, which satisfies (1) and (2), and tilt it near $D_i$ so
that the resulting $X$ becomes positively transverse to
$\mathcal{B}$, while still keeping properties (1) and (2).  (The
other option is to keep ${\bdry\over \bdry t}$ and smooth the spine
into a branched surface so that each $D_i$ no longer has any
vertical tangencies.) The foliation $\mathcal{F}$ is constructed by
first taking a lamination which is fully carried by $\mathcal{B}$
and by extending it to complementary regions which are $I$-bundles.
The $I$-fibers can be taken to be tangent to $X$ and hence the
foliations on the $I$-bundles transverse to $X$.
\end{proof}

Recall that the ambient manifold $M$ can be written as $N\cup
(S^1\times D^2)$, where the meridian of the solid torus has slope
$-{1\over k}$  on $\bdry N$, where $k$ is the integer in the proof
of Theorem~\ref{Roberts}.  The foliation $\mathcal{F}$ on $N$ is now
extended to all of $M$ (also called $\mathcal{F}$) by foliating
$S^1\times D^2$ by meridian disks.

\begin{lemma}
There exists an isotopy $\phi: S\times[0,1]\rightarrow N$ so that
the following hold:
\be
\item $\phi_s(S)$ is properly embedded for all $s\in[0,1]$.
Here $\phi_s(y) \stackrel{def} =\phi(y,s)$.
\item $\phi_0(S)=S\times\{0\}$.
\item $\phi_s(S)$ is positively transverse to $X$ for all $s$.
\item $\bdry (\phi_1(S))$ is positively transverse to $\mathcal{F}$.
\ee
\end{lemma}

Here the orientation on $\bdry (\phi_s(S))$ is the one induced from
$\phi_s(S)$, which in turn is consistent with that of
$S\times\{t\}$.

\begin{proof}
It suffices to show that there is an isotopy $\psi:
S^1\times[0,1]\rightarrow \bdry N$ so that:
\be
\item $\psi_s$ is an embedding for all $s\in[0,1]$.
\item $\psi_0(S^1)=\bdry (S\times\{0\})$ (and their orientations agree).
\item $\psi_s(S^1)$ is transverse to $X|_{\bdry N}$ for all $s$ and
$(\dot\psi_s,X)$ form an oriented basis of $\bdry N$.
\item $\psi_1(S^1)$ is positively transverse to $\mathcal{F}$. Here
$\mathcal{F}$ intersects $\bdry N$ transversely.
\ee

To demonstrate the existence of such an isotopy, we examine the
train track $\mathcal{T}=\mathcal{B}\cap \bdry N$, where
$\mathcal{B}$ is the branched surface constructed in Lemma~\ref{vf}.
We use standard Euclidean coordinates $(x,y)$ on $\bdry N\simeq
\R^2/\Z^2$ (given by Roberts as $(\lambda, \mu)$). By construction,
$\{y=0\}\subset \mathcal{T}$.  Now $\{-\varepsilon\leq y\leq
\varepsilon\}\cap \mathcal{T}$ has four branch points, all on
$\{y=0\}$.  Two of the branches come in from $\{0\leq y\leq
\varepsilon\}$ and the other two come in from $\{-\varepsilon\leq
y\leq 0\}$.  (See Figure~\ref{traintrack}.)  Since they are coming
from a single disk in $N$, the branching directions of the two
branches on $\{0\leq y\leq \varepsilon\}$ are opposite and so are
the branching directions of the two branches on $\{-\varepsilon\leq
y\leq 0\}$. Here the {\em branching direction} at a branch point is
the direction in which two branches come together to become one.  We
will assume that $X={\bdry \over \bdry y}$ on $\{-\varepsilon\leq
y\leq \varepsilon\}$.

\begin{figure}[ht]
\s
\begin{overpic}[width=13cm]{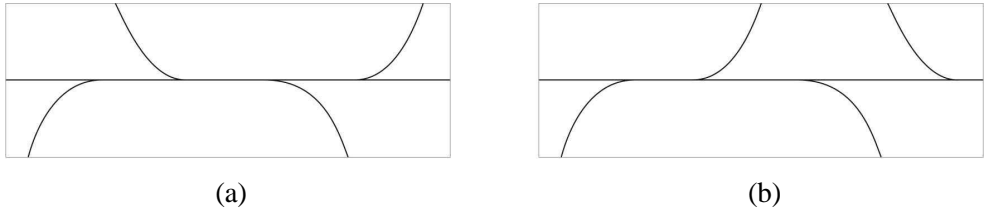}
\put(21.5,0){\small(a)}
\put(76,0){\small(b)}
\end{overpic}
\caption{The train track $\mathcal{T}$ near $y=0$. The bounding box
represents $[0,1]\times [-\varepsilon,\varepsilon]\subset
\R^2/\Z^2$.  Two of the possible combinations are labeled (a) and
(b).} \label{traintrack}
\end{figure}

We then let $N(\mathcal{T})$ be the train track neighborhood of
$\mathcal{T}$.  By Roberts' construction, $N(\mathcal{T})$ fully
carries a sublamination of $\mathcal{F}\cap T(\bdry N)$. Without
loss of generality, this sublamination $C$ satisfies the following:
\begin{enumerate}
\item[(i)] $C$ is a finite union of closed curves of slope $-{1\over k}$,
where $k$ is a positive integer.
\item[(ii)] The horizontal boundary of $N(\mathcal{T})$ is contained in $C$.
\end{enumerate}

Recall that $\lambda$ is oriented as $\bdry S$, and is directed by
${\bdry\over \bdry x}$.  Orient $\mu$ so that $\dot\mu$ has positive
${\bdry\over \bdry y}$--component.  Orient $\mathcal{T}$ (and hence
$C$) using the transverse vector field $X$.  More precisely, $(\dot
C,X)$ are to form an oriented basis for $\bdry N$.

\begin{figure}[ht]
\s
\begin{overpic}[width=13cm]{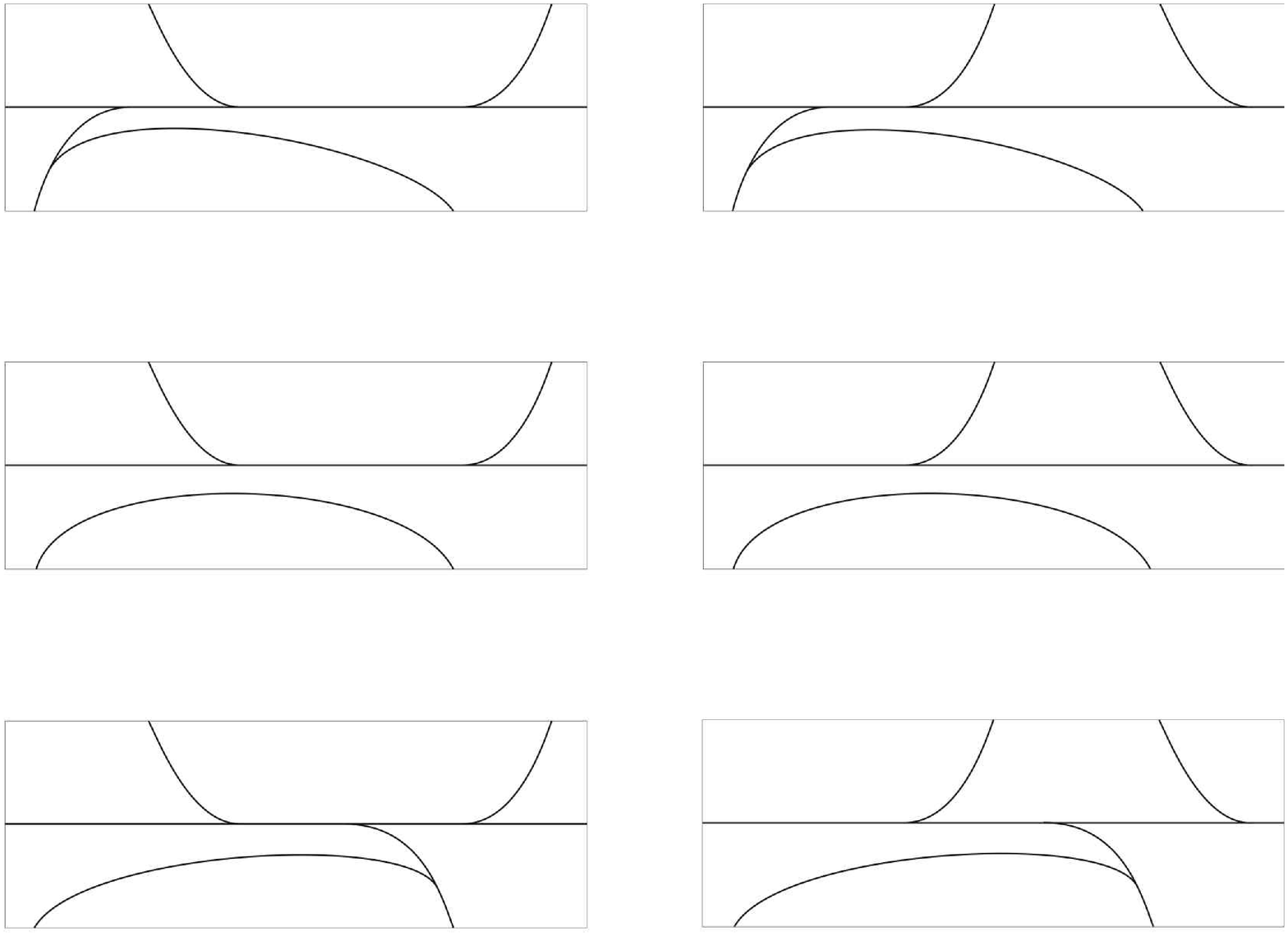}
\put(21.5,-1){\small($\mbox{a}_3$)}
\put(76,-1){\small($\mbox{b}_3$)}
\put(21.5,26.5){\small($\mbox{a}_2$)}
\put(76,26.5){\small($\mbox{b}_2$)}
\put(21.5,54.5){\small($\mbox{a}_1$)}
\put(76,54.5){\small($\mbox{b}_1$)}
\end{overpic}
\caption{Possible splittings of $\mathcal{T}$.}
\label{splitting1}
\end{figure}

We now split $\mathcal{T}$ by pushing in one of the branches in
$[0,1]\times [-\varepsilon,0]$ to obtain $\mathcal{T}_1$, which also
fully carries $C$ and satisfies (i) and (ii) above.  The
possibilities are given in Figure~\ref{splitting1}. We claim that
($\mbox{a}_1$), ($\mbox{a}_2$), ($\mbox{b}_1$), or ($\mbox{b}_2$)
are not possible for $\mathcal{T}_1$.  Indeed, in ($\mbox{a}_2$) and
($\mbox{b}_2$) the algebraic intersection number $\langle \lambda,
C\rangle=0$, so $\mbox{slope}(C)=0$, a contradiction.  In case
($\mbox{a}_1$) and ($\mbox{b}_1$), we have $\langle
\lambda,C\rangle>0$ and $\langle \mu,C\rangle<0$, implying that
$\mbox{slope}(C)>0$, which is also a contradiction.  Therefore,
$\mathcal{T}_1$ must be ($\mbox{a}_3$) or ($\mbox{b}_3$).  Since
$\mathcal{T}_1$ fully carries $C$, in either case there must exist a
subarc $\delta:[0,1]\rightarrow \R^2/\Z^2$ of the horizontal
boundary of $N(\mathcal{T}_1)$  such that $\delta(0)$ and
$\delta(1)$ have the same $x$-coordinate, $\delta$ ``winds around''
in the (positive) $x$-direction once, and the $y$-coordinate of
$\delta(1)$ is smaller than that of $\delta(0)$ (here we are in
$[0,1]\times[-\varepsilon,\varepsilon]$).  Let $\delta_1$ be the
oriented integral subarc of $X={\bdry\over \bdry y}$ in $[0,1]\times
[-\varepsilon,\varepsilon]$ from $\delta(1)$ to $\delta(0)$.  The
concatenation $\delta*\delta_1$ can easily be perturbed into a
closed curve which is isotopic to $\lambda=\bdry (S\times\{0\})$ and
is positively transverse to $\mathcal{F}$.  Moreover, it is easy to
take the isotopy to be transverse to $X$ throughout.
\end{proof}

Next, following Eliashberg and Thurston \cite{ET}, take a
$C^0$-small perturbation of $\mathcal{F}$, which we denote by
$\xi_{\mathcal{F}}$.  The characteristic foliation of
$\xi_{\mathcal{F}}$ on $\bdry N$ can be taken to have slope
$-{1\over k}+\varepsilon$, where $\varepsilon$ is an arbitrarily
small positive number. We can choose the  perturbation so that the
characteristic foliation is nonsingular Morse-Smale with two closed
orbits, one attracting and one repelling, and   $-{1\over
k}+\varepsilon$ is the slope of the closed orbits. (Hence $\bdry N$
is a convex surface with two dividing curves of slope $-{1\over
k}+\varepsilon$.) We make the perturbation $\xi_{\mathcal{F}}$
sufficiently close to $\mathcal{F}$ so that $\bdry(\phi_1(S))$ is
positively transverse to $\xi_{\mathcal{F}}$.

Since $X$ is positively transverse to $\xi_{\mathcal{F}}$ and also
to $S'\stackrel{def}=\phi_1(S)$, it follows that the characteristic
foliation of $\xi_{\mathcal{F}}$ on $S'$ does not have any negative
singular points.  Therefore,
\begin{equation} \label{selflinking}
l(\bdry S')=  -e^+ + h^+  + e^- - h^- =  -e^+ + h^+ =- \chi(S') = 2g(S')-1.
\end{equation}
where $l(\bdry S')$ denotes the {\em self-linking number} of the
transverse knot $\bdry S'$ {\em with respect to $S'$}, $e^{\pm}$ and
$h^{\pm }$ are the numbers of positive and negative elliptic and
hyperbolic tangencies of the contact structure on $S'$, and the
genus $g(\Sigma)$ of a compact surface $\Sigma$ with boundary is the
genus of the closed surface obtained by capping off all the boundary
components with disks. A good reference for invariants of transverse
and Legendrian knots is \cite{Et}.

\begin{rmk}
The fact that $X$ is positively transverse to both $S'$ and
$\xi_{\mathcal{F}}$ does not imply that the dividing set
$\Gamma_{S'}$ is empty.  Although there are no negative singular
points in the characteristic foliation, closed orbits of Morse-Smale
type can function as sinks.  Hence we can have annular regions of
$S'_-$, where $S'\setminus \Gamma_{S'}=S'_+\cup S'_-$.
\end{rmk}

We next explain how to pass from $S'$ to a convex surface $S''$ with
Legendrian boundary.

\begin{lemma}\label{leg-boundary}
There exists a convex surface $T$ isotopic to $\bdry N$ inside an
$I$-invariant neighborhood of $\bdry N$ so that
$\Gamma_T=\Gamma_{\bdry N}$, and a convex surface $S''$ with
Legendrian boundary (isotopic to $S'$) so that $\bdry S''\subset T$
and
$$tb(\bdry S'')-r(\bdry S'') = l(\bdry S')=2g(S')-1.$$
\end{lemma}

Here $tb(\bdry S'')$ and $r(\bdry S'')$ are the {\em
Thurston-Bennequin invariant} and {\em rotation number} of the
Legendrian knot $\bdry S''$ {\em with respect to  $\bdry S''$}.

\begin{proof}
We isotop $\bdry N$ inside its invariant neighborhood to obtain the
convex surface $T$. More explicitly, we tilt $\bdry N$ near the
closed orbits of the nonsingular Morse-Smale characteristic
foliation so that $T$ has Legendrian divides in place of closed
orbits. Then isotop $S'$ to $S''$ which has Legendrian boundary and
such that $\bdry S''\subset T$. It is not hard to see that the
positive transverse push-off of $\bdry S''$ is transversely isotopic
to $\bdry S'$. Finally recall that if $\gamma^+$ is a positive
transverse push-off of $\gamma$, then
$tb(\gamma)-r(\gamma)=l(\gamma^+)$ (note the sign in front of
$r(\gamma)$ is negative, not positive).
\end{proof}

Now recall that
$$r(\bdry S'')= \chi(S''_+)-\chi(S''_-),$$
where $S_+''$ and $S_-''$ are the positive and negative regions of
$S''-\Gamma_{S''}$. By comparison with Lemma~\ref{leg-boundary},
which states that:
$$r(\bdry S'')= 1-2g(S'') + tb(\bdry S''),$$
we have $\chi(S''_+)=1-2g(S'')$ and $-\chi(S''_-)=tb(\bdry S'')$.
This implies the following:

\begin{cor}
$\Gamma_{S''}$ consists of $\bdry$-parallel dividing arcs and
curves, together with pairs of parallel essential closed curves.
\end{cor}

Note here that by a {\em $\bdry$-parallel dividing arc} we mean a
properly embedded arc that cuts off a disk whose interior intersects
no other components of the dividing set. In particular, the disks
cut off by $\bdry$-parallel dividing arcs are disjoint. {\em A
$\bdry$-parallel closed dividing curve} is a closed curve parallel
to the boundary.

Let us now rename $\bdry N$ and $S=S\times\{1\}$ so that the
following hold:
\begin{itemize}
\item[(i)] $\bdry N$ is convex, $\#\Gamma_{\bdry N}=2$,
and $\mbox{slope}(\Gamma_{\bdry N})=-{1\over k}+\varepsilon$, where
$\varepsilon$ is a small positive number.
\item[(ii)] $S=S\times\{1\}$ has boundary on $\bdry N$, and
$\Gamma_S$ consists of $\bdry$-parallel dividing arcs and closed
curves, together with pairs of parallel essential closed curves.
\item[(iii)] The solid torus $S^1\times D^2=M-N$ is the standard
neighborhood of a Legendrian curve.
\end{itemize}
We will now normalize $\Gamma_S$ in a manner similar to Section~7 of
\cite{HKM1}.

\begin{prop} \label{prop: only bdry parallel arcs}
There exists a convex surface isotopic to $S$ whose dividing set
only consists of $\bdry$-parallel arcs.
\end{prop}

\begin{proof}
Consider the cut-open manifold $S\times[0,1]$.  Here
$\Gamma_{S\times\{1\}}=\Gamma_S$ and
$\Gamma_{S\times\{0\}}=\psi(\Gamma_S)$. Since the monodromy map
$\psi$ is pseudo-Anosov, $\Gamma_S\not=\psi(\Gamma_S)$ unless
$\Gamma_S$ is a union of $\bdry$-parallel arcs and $\bdry$-parallel
closed curves.

We will first reduce to the case of such a union. If
$\Gamma_S\not=\psi(\Gamma_S)$, then, by Proposition 7.1 of
\cite{HKM1}, there exists a closed curve $\gamma$, possibly
separating, which intersects $\Gamma_{S\times\{i\}}$, $i=0,1$,
efficiently and such that $\#(\gamma\cap\Gamma_{S\times\{1\}})
\not=\#(\gamma\cap \Gamma_{S\times\{0\}})$. Now apply the Legendrian
Realization Principle to make $\gamma\times\{0,1\}$ Legendrian, and
apply the Flexibility Theorem to make $\gamma\times[0,1]$ convex
with Legendrian boundary.  By the Imbalance Principle of \cite{H1},
there must exist a bypass along $\gamma\times\{0\}$, say.  Let
$\mathcal{B}_\alpha$ be the bypass and $\alpha$ the arc of
attachment for the bypass.

Note that the condition that $\gamma$ intersect
$\Gamma_{S\times\{i\}}$ efficiently eliminates the possibility of a
trivial bypass. Hence we have the following possibilities:

\s\n (i) If $\alpha$ intersects three distinct dividing curves, then
attaching $\mathcal{B}_\alpha$ yields a convex surface $S'$ isotopic
to $S$ with fewer dividing curves.

\s\n (ii) If $\alpha$ starts on a dividing curve $\gamma_1$, passes
through a parallel dividing curve $\gamma_2$, and ends on
$\gamma_1$, then $\gamma_1$ and $\gamma_2$ are nonseparating, and we
may apply Bypass Rotation (see \cite{HKM2}) so that one of the
endpoints of $\alpha$ ends on a different dividing curve $\gamma_3$
(here $\gamma_3$ may be a $\bdry$-parallel arc).  Then apply case
(i).

\s\n (iii) Suppose $\alpha$ starts on $\gamma_1$, passes through a
parallel $\gamma_2$, and ends on $\gamma_2$ after going around a
nontrivial loop.  There are two possibilities: either $\gamma_1$ and
$\gamma_2$ are both separating curves or they are both nonseparating
curves. If $\gamma_1$ and $\gamma_2$ are both nonseparating, then we
can apply Bypass Rotation and get to (i) and reduce the number of
dividing curves. If $\gamma_1$ and $\gamma_2$ are both separating,
but the connected component of $S\setminus \gamma_2$ containing the
subarc of $\alpha$ from $\gamma_2$ to itself has other components of
$\Gamma_S$, then we can apply Bypass Rotation, get  to (i) and
reduce the number of dividing curves.  Finally, if $\gamma_1$ and
$\gamma_2$ are separating and $\gamma_2$ splits off a subsurface of
$S$ which does not contain other components of $\Gamma_S$, then
attaching $\mathcal{B}_\alpha$ yields a pair of parallel dividing
curves which are either nonseparating or are separating but split
off a strictly smaller subsurface.  Hence, we can reduce the
complexity in one of two ways: either by reducing number of
separating curves or by reducing the genus of the separated part. By
repeating this procedure, we can reduce the genus of the separated
part down to $1$ in finitely many steps and force the appearance of
a nonseparating pair. Then apply case (i) or (ii). \s

We can repeat this procedure until we eliminate all pairs of curves
that are not parallel to the boundary.

To eliminate the closed curves parallel to $\bdry S$, we cut $N$
open along $S$ to obtain $S\times[0,1]$. The dividing set on $\bdry
(S\times[0,1])$, after rounding, will consist of $4n+1$ closed
curves which are parallel to $\bdry S$. There are $2n$ each on
$S\times\{0\}$ and $S\times\{1\}$, and one which is created from the
$\bdry$-parallel arcs by edge-rounding. Let us number the dividing
curves consecutively (as they appear on $\bdry (S\times[0,1])$) as
$\gamma_1,\dots,\gamma_{4n+1}$. Now let $\delta$ be a properly
embedded, non-boundary-parallel arc from $\bdry S$ to itself.  Then
cut $S\times[0,1]$ along the disk $D=\delta\times[0,1]$, which we
take to be convex with Legendrian boundary.  Furthermore, we take
$\bdry D$ to be efficient with respect to $\Gamma_{\bdry
(S\times[0,1])}$.  Now consider the $\bdry$-parallel arcs of
$\Gamma_D$.  The only time a $\bdry$-parallel arc does not have a
corresponding bypass which reduces $\#\Gamma_{S\times\{0\}}$ or
$\#\Gamma_{S\times\{1\}}$ or puts us in case (iii) above is if it
straddled the middle curve $\gamma_{2n+1}$. In this case, there are
only two $\bdry$-parallel arcs on $D$ and all the other dividing
arcs on $D$ are ``parallel'' to these $\bdry$-parallel arcs that
straddle $\gamma_{2n+1}$.

However, we claim that this particular form of $\Gamma_D$ implies
that there is a closed Legendrian curve $\delta$ which is isotopic
to a meridian curve on $\bdry N$ and has zero relative
Thurston-Bennequin invariant with respect to the tangent framing of
$\bdry N$. In fact, any properly embedded Legendrian arc $\delta_i$
on $D$ which is parallel to and disjoint from arcs of $\Gamma_D$,
and has endpoints on $\gamma_{i}$ and $\gamma_{4n+2-i}$, glues to
give such a closed Legendrian curve $\delta$, after possibly sliding
an endpoint along $\gamma_i$.  Now, $\delta$ bounds an overtwisted
disk in $M$, obtained from $N$ by Dehn filling along the meridian
slope. This contradicts the fact that the contact structure
$\xi_{\mathcal{F}}$ on $M$ is a tight contact structure. (Recall
that $\xi_{\mathcal{F}}$ is a perturbation of a taut foliation
$\mathcal{F}$.) This proves that $\Gamma_S$, after successive bypass
attachments, can be made to consist only of $\bdry$-parallel arcs.
\end{proof}

Finally, to prove that $\xi_{\mathcal{F}}$ is the same as the
contact structure defined by $(S,h)$, we cut $N$ open along $S$ with
only $\bdry$-parallel dividing arcs, as furnished by
Proposition~\ref{prop: only bdry parallel arcs}, and consider the
contact structure on the cut-open manifold induced by
$\xi_{\mathcal{F}}$. The dividing set for $\xi_{\mathcal{F}}$ on
$\bdry (S\times[0,1])$ consists of one closed curve parallel to
$\bdry S$. We can now take a system of arcs $\alpha_i$,
$i=1,\dots,2g(S)$, on $S$ and cut along the disks $\alpha_i
\times[0,1]$ to decompose $S\times[0,1]$ into a disk times $[0,1]$.
Note that each decomposition is along a convex disk $D_i$ isotopic
to $\alpha_i\times[0,1]$ with Legendrian boundary, so that $\bdry
D_i$ intersects $\Gamma_{\bdry (S\times[0,1])}$ in exactly two
points. Hence the dividing set $\Gamma_{D_i}$ is determined, i.e.,
consists of a single dividing arc. There is a unique (up to isotopy)
contact structure with such a decomposition, called a {\em product
disk decomposition}. The decomposition of $M$ into such an $N$ and a
standard neighborhood $S^1\times D^2$ of a Legendrian curve is
clearly identical to the contact structure adapted to the open book
$(S,h)$.

\s\n {\em Acknowledgements.} We thank Francis Bonahon and Bill
Thurston for helpful discussions. Discussions with \'Etienne Ghys
also helped clarify aspects of Section~\ref{braidgroup}.  We are
also very grateful to the referee for the extensive list of
comments.

\end{document}